\input amstex
\documentstyle{amsppt}
\vsize= 52pc
\NoBlackBoxes

\topmatter 
\title Quasiconvexity and amalgams\endtitle
\rightheadtext {Quasiconvexity and Amalgams}
\author Ilya Kapovich
\endauthor
\address Department of Mathematics, Graduate School and University Center of the City University of New York, 33 West 42-nd Street, New York, NY10036\endaddress
\email ilya\@groups.sci.ccny.cuny.edu\endemail
\subjclass Primary 2.2F10; Secondary 2.2F32 \endsubjclass
\abstract We obtain a criterion for quasiconvexity of a subgroup of an amalgamated free product of two word hyperbolic groups along a virtually cyclic subgroup. The result provides a method of constructing new word hyperbolic group in class (Q), that is such that all their finitely generated subgroups are quasiconvex. It is known that free groups, hyperbolic surface groups and most 3-dimensional Kleinian groups have property (Q). We also give some applications of our results to one-relator groups and exponential groups.
    
\endabstract
\thanks The author is supported by Alfred P. Sloan Foundation Doctoral Dissertation Fellowship
\endthanks

\endtopmatter
\hsize=6.5 true in

\head 0.Introduction\endhead

An important notion in the theory of word hyperbolic groups is the concept of a quasiconvex (or rational) subgroup, which, roughly speaking, corresponds to a geometrically finite subgroup of a classical hyperbolic group.
\definition {Proposition-Definition} (see [1] for proof)
Let $G$ be a word hyperbolic group and $A$ be a subgroup of $G$. Then the following conditions are equivalent:
\roster
\item "(1)" for some finite generating set $X=S\cup S^{-1}$ of $G$ there is an $\epsilon \ge 0$ such that for any $a\in A$ for each $d_X$-geodesic path $p$ from $1$ to $a$ in the Cayley graph $\Gamma (G,X)$ of $G$ for any point $x$ on $p$ there is $a^{\prime}\in A$ such that $d_X(x,a^{\prime})\le \epsilon$ (here $d_X$ denotes the word metric on $\Gamma (G,X)$ corresponding to $X$);
\item "(2)" $A$ is finitely generated and for some finite generating set $Y=T\cup T^{-1}$ of $A$ and for some finite generating set $X=S\cup S^{-1}$ of $G$ there is a constant $C>0$ such that for any $a\in A$
$$  d_Y(a,1)\le C\cdot d_X(a,1)+C.$$ 
\endroster
If any of these conditions is satisfied then $A$ is called a {\it quasicovex subgroup} of $G$.
\enddefinition
It can be shown that if $A$ is a quasiconvex subgroup of $G$ then conditions (1) and (2) of the previous definition are satisfied for any finite generating set of $G$ and any finite generating set of $A$.
Quasiconvex subgroups of word hyperbolic groups are themselves word hyperbolic and an intersection of a finite number of quasiconvex subgroups is again quasiconvex. Also finite subgroups, subgroups of finite index, virtually cyclic subgroups, free factors and conjugates of quasiconvex subgroups of word hyperbolic groups are quasiconvex (see [1], [8], [9], [10]).

The following class of groups is of considerable interest.
\definition{Definition} We say that a word hyperbolic group $G$ has {\it property (Q)} if any finitely generated subgroup of $G$ is quasiconvex in $G$.
\enddefinition

We note some good properties of groups with property (Q).
\roster
\item a finitely generated subgroup of a group with property (Q) also has property (Q) and so it is finitely presented and word hyperbolic;
\item a finite extension of a group with property (Q) also has property (Q);
\item if $G$ has property (Q) then the intersection of any two finitely generated subgroups of $G$ is finitely generated (that is $G$ has the {\it Howson property}) (see [19]);
\item  any infinite finitely generated subgroup $H$ of $G$ has finite index in its virtual normalizer $VN_G(H)=\{g\in G\mid |H:gHg^{-1}\cap H|<\infty, |gHg^{-1}:gHg^{-1}\cap H|<\infty\}$ (see [14]);
\item if $H_1,H_2$ are infinite subgroups of $G$ such that $H=H_1\cap H_2$ has finite index in both $H_1$ and $H_2$ then $H$ is of finite index in $E=gp(H_1,H_2)$ (see [14]);
\item if $G$ is a torsion-free group from class (Q) then any maximal cyclic subgroup of $G$ is a {\it Burns subgroup} (see [4] and [6] for definitions) of $G$ (see [12] for proof). In particular it implies that 
if $B$ is a group with the Howson property, $b\in B$ is an element of infinite order, $g\in G$ is an element of infinite order which is not a proper power then the group $G\underset{g=b}\to\ast B$ has the Howson property.
 
\endroster

It seems that most word hyperbolic groups have property (Q) but, nevertheless, there are relatively few examples for which it is proven. A finitely generated free group and  the fundamental group of a closed hyperbolic surface belong to class (Q) (see [19], [18] and [20]). Also, if $G$ is a torsion-free geometrically finite Kleinian group without parabolics whose limit set is not the whole $S^2$, then $G$ has property (Q) (see [20]). 

In this paper we show how to build new groups with property (Q) from existing ones using free constructions.

\proclaim {Theorem A} 
Suppose $G=A_{-1}\ast_C A_1$ is word hyperbolic group where $C$ is virtually cyclic and the groups $A_1$, $A_{-1}$ have property (Q).
Then $G$ has property (Q).
\endproclaim

In order to obtain this result we need the following statement which characterizes quasiconvex subgroups of an amalgamated free product of two word hyperbolic groups and is of considerable independent interest.

\proclaim {Theorem B} 
Let $G=A_{-1}\ast_C A_1$ be a word hyperbolic group where $C$ is virtually cyclic (this implies, by the results of [11], that $A_1$ and $A_{-1}$ are quasiconvex in $G$).
Suppose $H$ is a finitely generated subgroup of $G$. 

Then the following conditions are equivalent.
\roster
\item "(a)" $H$ is quasiconvex in $G$;
\item "(b)" for any $g\in G$ and for any $i=1,-1$ the subgroup $gHg^{-1}\cap A_i$ is quasiconvex in $A_i$.
\endroster
\endproclaim 

\proclaim {Corollary 1} Let $A_1$, $A_{-1}$ be groups with property (Q) and suppose that $A_1$ is torsion-free. Let $x\in A_1$ be an element of infinite order which is not a proper power. Let $y\in A_{-1}$ be an element of infinite order.
Then the group $A_1\underset{x=y}\to\ast A_{-1}$ is word hyperbolic and has property (Q).
\endproclaim

\proclaim {Corollary 2} Let $G=A_1\ast_C A_{-1}$ where $C$ is finite and $A_1$, $A_{-1}$ have property (Q). Then $G$ is word hyperbolic and also has property (Q).
\endproclaim
Corollary 2 is a generalization of the result of [11] where its statement was proved when $C=\{1\}$.

There are also some interesting consequences of these results for so-called {\it exponential} groups (see [15] for definitions).

\proclaim {Corollary 3} Let $G$ be a torsion-free hyperbolic group with property (Q) (e.g. finitely generated free group, hyperbolic surface group etc).
Let $G^{\Bbb Q}$ be the tensor ${\Bbb Q}$-completion of $G$ where ${\Bbb Q}$ is the ring of rational numbers.
Then
\roster
\item  $G^{\Bbb Q}$ is a locally (Q)-group that is any finitely generated subgroup of $G^{\Bbb Q}$ is word hyperbolic and has property (Q);
\item $G^{\Bbb Q}$ has the Howson property that is the intersection of any two finitely generated subgroups of $G$ is finitely generated
\item if $H_1$ and $H_2$ are infinite commensurable subgroups of $G$, that is the intersection $H=H_1\cap H_2$ has finite index in both $H_1$ and $H_2$, then $H$ has finite index in their join $E=gp(H_1\cup H_2)$.
\endroster
\endproclaim

As a by-product of our results we also obtain the following statement.

\proclaim {Corollary 4} Let $G=G_1\ast_C G_2$ be a word hyperbolic group such that $C$ is finitely generated. Suppose $H$ is a finitely generated subgroup of $G$ such that for any $g\in G$ we have $g^{-1}Hg\cap G_1=g^{-1}Hg\cap G_2=\{1\}$. Then $H$ is quasiconvex in $G$.
\endproclaim

\proclaim {Corollary 5} Suppose $G$ is a one-relator group $G=<x_1,\dots ,x_k,y_1,\dots ,y_s |vu=1>$ where $v$ is a nontrivial freely reduced word in $x_1,\dots ,x_k$, $u$ is a nontrivial freely reduced word in $y_1,\dots ,y_s$ and $u$ is not a proper power in the free group $F(y_1,\dots ,y_s)$. Then $G$ has property (Q).
\endproclaim
Notice that fundamental groups of closed hyperbolic surfaces have one-relator presentations as in Corollary 5. G.Swarup [20] and C.Pittet [18] showed using the techniques of hyperbolic geometry that these groups belong to class (Q). Corollary 5 gives another, more combinatorial, proof of this fact.

It follows from the result of R.Burns [4] that a group $G$ as in Corollary 5 has the Howson property. Since the groups from class (Q) have the Howson property, Corollary 5 may be considered as a generalization of Burns' theorem. 

We would like to stress that quasiconvexity of a subgroup is not  a question of the isomorphism type of the subgroup but rather that of comparing the word metrics on the subgroup and on the ambient group. This makes the proof of Theorem B rather more difficult than it may seem from the first sight. 
To illustrate this point, consider the following example. Let $M$ be a closed hyperbolic 3-manifold fibering over a circle with fiber $S$, where $S$ is a closed hyperbolic surface. We may also think of $G$ as a geometrically finite group of isometries of ${\Bbb H}^3$ such that ${\Bbb H}^3/G=M$.
Let $x_0\in S$ and $G=\pi_1(G,x_0)$, $H=\pi_1(S,x_0)$.
Then there is a short exact sequence
$$1\rightarrow H{\rightarrow} G\rightarrow {\Bbb Z}\rightarrow 1$$
and therefore $H$ is not quasiconvex in $G$ (see [1]). Consider a simple closed curve $\gamma$ on $S$ passing through $x_0$ such that $\gamma$ separates $S$ into two non-contractible components. Then $\gamma$ defines a decomposition of $H$ as an amalgamated free product $H=F_1\ast_C F_2$ where $F_1, F_2$ are nonabelian free groups, $C$ is an infinite cyclic group which is malnormal in both $F_1$ and $F_2$.
It follows from geometric considerations that both $F_1$ and $F_2$ are geometrically finite groups of isometries of ${\Bbb H}^3$ and therefore (see [20]) both $F_1$ and $F_2$ are quasiconvex in $G$.
Thus we see that $F_1, F_2$ and $C$ are quasiconvex in $G$ and $H=gp(F_1,F_2)\simeq F_1\ast_C F_2$ is not quasiconvex in $G$.

\head 1.Subgroup structure of an amalgamated product \endhead

\subhead Some definitions and notations \endsubhead

If $G$ is a finitely generated group and $X$ is a finite generating set of $G$ closed under taking inverses, we denote the {\it Cayley graph} of $G$ with respect to $X$ by $\Gamma(G,X)$. The word metric on $\Gamma(G,X)$ corresponding to $X$ is denoted $d_X$.
Also, for an element $g\in G$ we put $l_X(g)=d_X(g,1)$.
If $w$ is a word in $X$, we denote by $\overline w$ the element of $G$ represented by $w$.
A word $w$ in $X$ is termed $d_X$-{\it geodesic} if the length $l(w)$ of $w$ is equal to $l_X(\overline w)$.
A word $w$ in $X$ is called $\lambda$-{\it quasigeodesic} with respect to $d_X$ if for any subword $u$ of $w$ we have $l(u)\le \lambda\cdot l_X(\overline u) +\lambda$.

We also will need some notations regarding graphs of groups.
Let ${\Bbb A}$ be a graph of groups and $A$ be the underlying graph of ${\Bbb A}$.
Then $VA$ and $EA$ denote the set of vertices and the set of edges of $A$ respectively.
We also denote by $E^{+}(A)$ the set of positively oriented edges of $A$. If $e$ is an edge of $A$ then its inverse is denoted by $e^{-1}$.
 For any vertex $v$ of $A$ the corresponding vertex group is denoted $A_v$. Similarly , if $e$ is an oriented edge of $A$, the edge group corresponding to $e$ is denoted $A_e$.
We also denote the initial vertex of $e$ by $\partial_0(e)$ and the terminal vertex of $e$ by $\partial_1(e)$. The edge monomorphism $A_e\rightarrow A_{\partial_0(e)}$ is denoted by $\alpha_e$. The edge-monomorphism $A_e\rightarrow A_{\partial_1(e)}$ is denoted by $\omega_e$. Recall that $\partial_0(e)=\partial_1(e^{-1})$ and $(e^{-1})^{-1}=e$ for any $e\in EA$. We also have $A_e=A_{e^{-1}}$ and $\alpha_e=\omega_{e^{-1}}$ for every $e\in EA$.

\subhead The premises \endsubhead

Suppose $G$ is a finitely generated group and
$$G=A_1\underset {C}\to\ast A_{-1} \eqno (1)$$ where $A_1, A_{-1}$ and $C$ are finitely generated.

Let $X_{-1}$ and $X_1$ be finite generating sets for $A_{-1}$ and $A_{-1}$ closed under inversions and containing a finite generating set ${\Cal C}$ of $C$. 
 Put ${\Cal G}=X_{-1}\cup X_1$ to be a finite generating set for $G$.
We denote by $d_{X_i}$ the word metric corresponding to $X_i$ on $A_i$, $i=\pm 1$.
Denote by $d_{\Cal G}$ the word metric corresponding to ${\Cal G}$ on $G$.
Also, fix an ordering on the sets $X_1$, $X_{-1}$.

Let $L_i$ be the set of all $X_i$-geodesic words $w$ such that
\roster 
\item "(i)" $\overline w$ is shortest (with respect to $d_{X_i}$) in the coset $\overline wC$;
\item "(ii)" if $v$ is another $X_i$-geodesic word such that $\overline vC=\overline wC$
 and $l(v)=l(w)$ then $w$ is lexicographically smaller than $v$.
\endroster

Then $T_i=\overline {L_i}$ is a left transversal for $C$ in $A_i$.

An expression of $g$ as an alternating product $g=u_1\dots u_k$ where $u_j\in A_{-1}\cup A_1$, $u_j\not\in C$
for $j<k$, and $u_{s+1}\not\in A_i$ whenever $u_s\in A_i$, $s=1,\dots ,k-1$, 
is called a {\it reduced form} of $g$ with respect to presentation (1). The elements $u_i$ are called {\it syllables} of $g$.
If $g=u_1\dots u_k$ is a reduced form of $g$ and $u_k\not\in C$ then we say
that $g$ has syllable length $k$. If $u_k\in C$ and so $g=u_k$, we say that
$g$ has syllable length zero.

If $x=u_1\dots u_k$ and $y=v_1\dots v_s$ are reduced forms of $x$ and $y$,
we say that $y$ is a {\it right segment} of $x$ and that
$x$ {\it ends} in $Cy$ if $s\le k$ and $u_{k-s+1}\dots u_ky^{-1}\in C$.
Under these circumstances $y$ is said to be a {\it proper right segment} of $x$ if $s<k$ or $y\in C$ and $x\not\in C$.
It is not hard to see that these definitions do not depend on the choices of reduced forms for $x$
and $y$.

Analogously if $x=u_1\dots u_k$ and $y=v_1\dots v_s$ are reduced expressions,
we say that $y$ is a {\it left segment} of $x$ and $x$ begins in $yC$ if $s\le k$ and $y^{-1}u_1\dots u_s\in C$.
A left segment $y$ of $x$ is said to be {\it proper} if $s<k$ or $y\in C$, $x\not\in C$.
Again these definitions do not depend on the choices of reduced forms for $x$ and $y$.
Observe also that $y$ is a (proper) right segment of $x$ if and only if $y^{-1}$ is 
a (proper) left segment of $x^{-1}$.
If $g=u_1\dots u_k$ is a reduced expression and $u_k\in A_i-C$, the element
$g$ is said to {\it end} in $A_i$.
Elements of $C$ are said to end in $C$.

\subhead A little bit of Bass-Serre theory \endsubhead

Let $H$ be a finitely generated subgroup of $G$ which is not elliptic, that is $H$ is not conjugate to a subgroup of $A_i$.

Let $\hat T$ be a Bass-Serre tree associated with the free product decomposition
$G=A_{-1}\ast_C A_1$. The vertices of $\hat T$ are just coset classes $gA_i$ of $A_1$ and
$A_{-1}$ in $G$. 
There is a distinguished vertex $d_1=A_1$ which is a basepoint of $\hat T$. Also denote $d_{-1}=A_{-1}$.
For any vertex $gA_i$ and $a\in T_i-\{1\}$ there is a positively oriented edge $(gA_i, gaA_{-i})$ which we label by $a$. There is also a positively oriented edge $(d_1, d_{-1})=(A_1,A_{-1})$ labelled by $1$.
The action of $G$ on $\hat T$ is obvious: $g\cdot fA_i=gfA_i$, $g\cdot (fA_i,faA_{-i})=(gfA_i,gfaA_{-i})$.  
The stabilizer in $G$ of a vertex $gA_i$ is clearly $gA_ig^{-1}$
and the stabilizer in $G$ of an edge $(gA_i,ga_iA_{-i})$ is equal to
$gA_ig^{-1}\cap ga_iA_{-i}a_i^{-1}g^{-1}=ga_i(A_i\cap A_{-i})a_i^{-1}g^{-1}=ga_iCa_i^{-1}g^{-1}$.
We will say that a vertex $v=gA_i$ of $\hat T$ has {\bf type} $A_i$. Any edge-path $p=(e_1,\dots ,e_k)$ such that each $e_j$ is positively oriented and labelled by $a_j$, has a label $a_1\dots a_k\in G$. For any vertex $v$ of $T$ there is a unique reduced edge-path $p_v$ from $d_1$ to $v$ whose label is denoted by $s_v$. Notice that $s_v=1$ if and only if $v=d_{\pm 1}$. We say that every vertex $w$ on $p_v$ is {\it less than or equal to} $v$  and write $w\le v$. It is obvious that $"\le "$ is a partial ordering on $V\hat T$. For a vertex $v\ne d_1$ the closest to $v$ vertex on $p_v$ which is different from $v$ is called the {\it preceding vertex} for $v$. In other words, $u$ is a preceding vertex for $v$ if and only if $(u,v)$ is a positively oriented edge of $\hat T$.
\smallskip
Then $H$ acts on $\hat T$ as a subgroup of $G$ and two vertices $g_1A_i$ and $g_2A_i$
lie in the same $H$-orbit if and only if the double coset classes
$Hg_1A_i$ and $Hg_2A_{-i}$ are equal.
There is a subtree $T$ of $\hat T$ which is $H$-invariant and does not contain any proper $H$-invariant subtrees, that is the action of $H$ on $T$ is {\it minimal}.
If $T$ is a single vertex, say $T=gA_i$ then $H\le gA_ig^{-1}$ which is impossible since we assumed that $H$ is not conjugate to a subgroup of $A_i$.
Thus $T$ has at lease one edge.
There is a finite subtree $Y$ of $T$ which serves as a "fundamental domain" for the action of $H$, that is any edge of $T$ lies in the $H$-orbit of a unique edge of $Y$.
Then we can find a subtree $Y_1$ of $Y$ such that 
any vertex of $T$ is $H$-equivalent a unique vertex of $Y_1$.
 Thus $Y$ is a union of $Y_1$ and a finite number of disjoint edges.
By conjugating $H$ we may assume that the edge $(d_{-1},d_1)$ is in $Y_1$ where $d_{-1}=A_{-1}, d_1=A_1$.
Notice that for $v\in VY, v\ne d_1$ the preceding vertex for $v$ belongs to $Y_1$.

Observe also that no edge $(u,v)$ of $Y$, where $u$ precedes $v$, other than $(A_{1},A_{-1})$, has label $1$.
Clearly, if $v$ is of type $A_i$ then $v=s_vA_i$.
Now let $v=s_vA_i$ be a vertex of $Y$. Put $A_v=A_i\cap s_v^{-1}Hs_v$.
Thus $H\cap s_vA_is_v^{-1}=s_vA_vs_v^{-1}$.
Also, let $e=(s_uA_{-i},s_vA_i)=(s_uA_{-i},s_ua_{-i}A_i)$ be a positively oriented edge of $Y$.
Recall that its stabilizer in $G$ is $s_vCs_v^{-1}$.
Put $C_e=C\cap s_v^{-1}Hs_v$.
Then $s_vC_es_v^{-1}=s_vCs_v^{-1}\cap H=s_v(A_v\cap C)s_v^{-1}$ is a subgroup of $s_uA_us_u^{-1}$ and $s_vA_vs_v^{-1}$.
Notice also that if $v=s_vA_i$ is a vertex of $Y-Y_1$ and $q=s_qA_i$
is the only vertex of $Y_1$ $H$-equivalent to $v$ then for any $h\in H$ such that $hq=v$ we have $h^{-1}s_vA_vs_v^{-1}h=s_qA_qs_q^{-1}$.

Suppose now that $v=s_vA_i$ is a vertex of $Y-Y_1$ and $s_v=s_ua_{-i}$ where
$u$ is the vertex of $Y_1$ preceding $v$ and $a_{-i}\in T_{-i}-\{1\}$. Let $q=s_qA_i$ be the only vertex of $Y_1$ which is $H$-equivalent to $v$. Then $Hs_qA_i=Hs_vA_i$, so there is an element $a\in A_i$ such that $s_v^{-1}s_qa\in H$.
We fix this element $a\in A_i$ for each $v\in Y-Y_1$ and denote $h_v=s_va^{-1}s_q^{-1}\in H$.

Clearly $s_va^{-1}s_q^{-1}\cdot q=s_va^{-1}s_q^{-1}\cdot s_qA_i=s_vA_i=v$, that is $h_vq=v$.
Since $h_vq=v$ we have $h_v^{-1}s_vA_vs_v^{-1}h_v=s_qA_qs_q^{-1}$.
Then $s_v(A_v\cap C)s_v^{-1}$ is a subgroup of $s_uA_us_u^{-1}$ and 
$h_v^{-1}s_v(A_v\cap C)s_v^{-1}h_v=s_qas_v^{-1}s_v(A_v\cap C)s_v^{-1}s_va^{-1}s_q^{-1}=s_qa(A_v\cap C)a^{-1}s_q^{-1}\le s_qA_qs_q^{-1}=s_qAs_q^{-1}\cap H$.
Thus the element $h_v$ conjugates the subgroup $s_v(A_v\cap C)s_v^{-1}$ of $A_u$
into the subgroup $s_qa(A_v\cap C)a^{-1}s_q^{-1}$ of $s_qA_qs_q^{-1}$.

The quotient graph of groups for the action of $H$ on $T$ can obtained from $Y$ in the following way. Let $B$ be an oriented graph such that
\roster 
\item "(1)" the vertices of $B$ are the vertices of $Y_1$;
\item "(2)" $B$ has one edge $(v,w)$ for each positive edge $(v,w)$ of $Y_1$;
\item "(3)" for any vertex $v\in V(Y-Y_1)$ and a vertex $v\in VY_1$ which is $H$-equivalent to $v$ there is an oriented edge $(u,q)$ in $B$ where $u$ is a preceding vertex for $q$.
\endroster

We give $B$ the structure of a graph of groups in the following way.
For any $v\in VY_1$ put $B_v=s_vA_vs_v^{-1}$ to be the vertex group of $v$.
For each edge $e=(u,v)$ of $Y_1$ where $u$ precedes $v$, put the edge group $B_e=s_v(A_v\cap C)s_v^{-1}$ where corresponding edge homomorphisms $\alpha_e\colon s_v(A_v\cap C)s_v^{-1}\rightarrow s_vA_vs_v^{-1}$ and $\partial_1\colon s_v(A_v\cap C)s_v^{-1}\rightarrow s_uA_us_u^{-1}$ are just the inclusion maps.

For any $v\in V(Y-Y_1)$, $u\in VY_1$ preceding $v$ and $q\in VY_1$ which is $H$-equivalent to $v$, put $B_e=s_v(A_u\cap C)s_v^{-1}\le s_uA_us_u^{-1}$ where $e=(u,q)\in EB$.
The boundary homomorphism $\alpha_e\colon B_e=s_v(A_u\cap C)s_v^{-1}\rightarrow s_uA_us_u^{-1}=B_u$ is the inclusion map.
The boundary homomorphism $\omega_e\colon B_e==s_v(A_u\cap C)s_v^{-1}\rightarrow s_qA_qs_q^{-1}=B_q$ is conjugation by $h_v$.
That is $\omega_e(g)=h_v^{-1}gh_v$ for any $g\in B_e$.
This defines a graph of groups ${\Bbb B}$.
Notice that $Y_1$ is a maximal subtree of $B$.
The fundamental group of the graph of groups ${\Bbb B}$ with respect to the maximal subtree $Y$ has the presentation

$$\pi_1({\Bbb B}, Y_1)=\underset{v\in VY_1}\to {(\ast B_v)}\ast F(E^{+}B)/\{e=1, e\in EY_1; \alpha_e(b)e=e\omega_e(b), e\in E^{+}B, b\in B_e\} \eqno (2) $$

Then by the fundamental result of Bass-Serre theory the map
$f\colon \pi_1({\Bbb B}, Y_1) \rightarrow H$ defined by
$f(g)=g$ for any $g\in B_v=s_vA_vs_v^{-1}$, $v\in VY_1$,
$f(e)=h_v$ where $e=(u,q)\in E^{+}(B-Y_1)$, $u$ precedes $v\in V(Y-Y_1)$, $u\in VY_1$ is $H$-equivalent to $v$,
is an isomorphism.
We will identify $\pi_1({\Bbb B}, Y_1)$ with $H$ via this isomorphism and will right
$$H=\pi_1({\Bbb B}, Y_1) \eqno (3)$$

\subhead Normal forms for the fundamental group of a graph of groups
\endsubhead

Let ${\Bbb A}$ be the graph of groups with underlying graph $A$ and let $T_0$ be the maximal subtree of $A$. Let $d_0$ be a fixed vertex of $T_0$. We will describe a set of normal forms for the fundamental group of ${\Bbb A}$

$$G_0=\pi_1({\Bbb A}, T_0) \eqno (4)$$  
which is slightly non-standard but which is more suitable for our purposes.

\definition{Definition 1.1}
A sequence

$$p=(g_1,e_1,g_2,e_2,\dots ,g_k,e_k,g_{k+1})$$

is called a {\it reduced sequence} if
\roster
\item $e_1\dots e_k$ is an edge-path in $A$;
\item for each $i=1,\dots k$ $g_i\in A_{\partial_0(e_i)}$ and
$g_{i+1}\in A_{\partial_1(e_i)}$;
\item $\partial_0(e_1)=\partial_1(e_k)=d_0$;
\item $e,1,e^{-1}$ is not a subsequence of $p$;
\item for any $i=2,\dots ,k+1$ either $g_i=1$ or $g_i\not\in \omega_{e_{i-1}}(B_{e_i})$;
\item if $i<j$, $g_i\ne 1$, $g_j\ne 1$,  $g_s=1$ for $i<s<j$  then

$$g_i\not\in 
           (\alpha_{i}\omega_{i}^{-1}) \dots (\alpha_{j-1}\omega_{j-1}^{-1})(B_{\partial_1(e_{j-1})})$$ 

\endroster

In the situation above we say that the number of terms in $p$ which are different from 1 in $G_0$ is
the {\it syllable length} of $p$. 
Any subsequence of $p$ represents an element of $G_0$, which is just the product of all terms in this subsequence viewed as elements of $G_0$. It is clear from the theory of graphs of groups that no subsequence of $p$ represents an element of a vertex group of $A$ unless this subsequence has  at most one term different from 1 in $G_0$.
If $g$ is the element of $G_0$ represented by $p$, we say that $p$ is a {\it reduced form} of $g$ with respect to presentation (4).
Let $p_1=(u_1,\dots ,u_n)$ be obtained from $p$ by deleting all terms which are equal to 1 in $G_0$. Thus each $u_i$ is either a stable letter or a nontrivial element of a vertex group. 
We call $p_1$ a {\it normal form} of $g$ with respect to presentation (4).
\enddefinition

\subhead Some Calculations\endsubhead

Recall that $H=\pi_1({\Bbb B}, Y_1)\le G=A_1\ast_C A_{-1}$.

\proclaim {Lemma 1.2}
\roster
\item "(1)" If $e=(s_vA_i,s_va_iA_{-i})=(s_vA_i,s_xA_{-i})$ is a positively oriented edge of $Y$, where $a_i\ne 1$,
then $a_iC\cap A_v=\emptyset$.
\item "(2)" If $v=s_vA_i$ is a vertex of $Y-Y_1$, $q$ is the only vertex of $Y_1$ which is $H$-equivalent to $v$ and $h_v=s_va^{-1}s_q^{-1}$
then $A_qa\cap C=\emptyset$.
\item "(3)" If $v_1,v_2$ are distinct vertices from $Y-Y_1$ which are $H$-equivalent to the same vertex $q$ of $Y_1$ and
$h_{v_j}=s_{v_j}a_j^{-1}s_q^{-1}$, $i=1,2$ then
$A_qa_1C\not=A_qa_2C$.
\item "(4)" If $q=s_qA_i\in Y_1$ is a preceding vertex for $w=s_qbA_{-i}\in VY$ and for some $v\in V(Y-Y_1)$ $h_v=s_va^{-1}s_q^{-1}$, $a\in A_i$ then $bC\cap A_qa=\emptyset$. 
\item "(5)" If $v=s_vA_{-i}=s_uaA_{-i}\in VY$, $w=s_wA_{-i}=s_ubA_{-i}\in VY$, $a\not=b, a,b\in T_i$, $u=s_uA_i$, where the edges $(u,v)$ and $(u,w)$ of $Y$ are positively oriented,  then $A_ua\not=A_ub$.

\endroster
\endproclaim

\demo {Proof}

(1) Suppose $a_ic\in A_v$, $c\in C$ that is $s_va_ics_v^{-1}=h\in H$.

Then $s_v=s_ua_{-i}$, where $u=s_uA_{-i}$ is a preceding vertex for $v$. Notice that $u,v\in Y_1$.
We have 
$$hu=s_va_ics_v^{-1}u=s_va_ica_{-i}^{-1}s_u^{-1}\cdot s_uA_{-i}=s_va_iA_{-i}=x.$$

On the other hand $hv=v$. Thus $h$ takes the edge $(u,v)$ into the edge $(x,v)$ what contradicts our assumptions that no two edges of $Y$ are $H$-equivalent.

(2) Let $v=s_vA_i\in Y-Y_1$ and let $q=s_qA_i$ be the vertex of $Y_1$ $H$-equivalent to $v$.
Let $s_v=s_ua_{-i}$, where $u=s_uA_{-i}$ is a vertex preceding to $v$. Suppose $a_qa=c\in C$ for some $a_q\in A_q$. Then $a^{-1}a_q^{-1}=c^{-1}$.
As before $h_v=s_va^{-1}s_q^{-1}\in H$ and $h_v\cdot q=v$.
Notice also that $h_0=s_va_q^{-1}s_v^{-1}\in H$, $h_0\cdot q=q$. Therefore $h\cdot q=v$ where $h=h_vh_0$.
Now let $y=s_yA_{-i}$ be the vertex preceding $q$ and $s_q=s_ya_{-i}^{\prime}$.
Then $h\cdot y=s_va^{-1}s_q^{-1}\cdot s_qa_q^{-1}s_q^{-1} \cdot s_yA_{-i}=s_ua_{-i}a^{-1}(a_{-i}^{\prime})^{-1}s_y^{-1}\cdot s_ya_{-i}^{\prime}a_q^{-1}(a_{-i}^{\prime})^{-1}s_y^{-1}\cdot s_yA_{-i}=s_v(a^{-1}a_q^{-1}){a_{-i}^{\prime}}^{-1}A_{-i}=s_vc^{-1}A_{-i}=s_vA_{-i}=s_ua_{-i}A_{-i}=s_uA_{-i}=u$.
Thus $h\cdot y=u$ and $h\cdot q=v$ and $h$ takes the edge $(y,q)$ into the edge $(u,v)$.
This contradicts our assumptions that no two edges of $Y$ are $H$-equivalent.

(3) Suppose $A_qa_1C=A_qa_2C$ that is $a_1^{-1}a_qa_2=c\in C$ for some $a_q\in A_q$. Put $h_0=s_qa_qs_q^{-1}\in H$.
Let $v_1=s_{v_1}A_i$, $v_2=s_{v_2}A_i$. Thus we know that $q=s_qA_i$ and $a_{\pm 1}\in A_i-C$.
Put $h=h_{v_1}h_0h_{v_2}^{-1}$.
Then $h\cdot v_2=v_1$ since $h_{v_j}(q)=v_j, j=\pm 1$ and $h_0(q)=q$.
Let $s_{v_j}=s_{u_j}b_j$ where $b_j\in T_{-j}$ and $u_j=s_{u_j}A_{-i}$ is the preceding vertex for $v_j$, $j=1,2$.

Then $hu_2=s_{u_1}b_1a_1^{-1}s_q^{-1}\cdot s_qa_qs_q^{-1}\cdot s_qa_2b_2^{-1}s_{u_2}^{-1}\cdot s_{u_2}A_{-i}=s_{u_1}b_1cb_2^{-1}A_{-i}=s_{u_1}A_{-i}=u_1$.
Thus $h$ takes the edge $(u_2,v_2)$ into the edge $(u_1,v_1)$ what contradicts our assumptions that no two distinct edges of $Y$ are $H$-equivalent.

(4) Suppose $bca^{-1}=a_q$, $a_q\in A_q, c\in C$.
Consider the preceding $v$ vertex $u=s_uA_{-i}\in VY_1$.
Then $s_v=s_uf$ for some $f\in T_{-i}$.
We have $h=s_qbca^{-1}s_q^{-1}\in H$ and $h_v^{-1}=s_qas_v^{-1}\in H$. Thus $h_1=hh_v^{-1}=s_qbcf^{-1}s_u^{-1}\in H$.
Clearly $h_1(u)=s_qbcf^{-1}s_u^{-1}s_uA_{-i}=s_qbA_{-i}=w$ and $h_1(v)=s_qbcf^{-1}s_u^{-1}s_ufA_i=s_qA_i=q$.
Thus $h_1$ takes the edge $(u,v)$ into the edge $(w,q)$ what contradicts the fact that no two distinct edges of $Y$ are $H$-equivalent.

(5) Suppose $a=a_ub$, $a_u\in A_u$. Then $h=s_uab^{-1}s_u^{-1}\in H$. However $h(u)=u$ and $h(w)=v$. Thus $h$ takes the edge $(u,w)$ into the edge $(u,v)$ which is impossible.

\enddemo

\proclaim {Lemma 1.3} If $v_1$, $v$ are vertices of $Y$ of type $A_i$ and $Hs_vC=Hs_{v_1}C$ 
then $v=v_1$. 
\endproclaim
\demo {Proof} Suppose $q=s_qA_i\ne v=s_vA_i$ and $Hs_qC=Hs_vC$, that is $s_qcs_v^{-1}=h$ for some $h\in H$, $c\in C$.
Let $us_uA_{-i}$ be the preceding vertex for $v$ when $v\ne d_1$ and let $u=d_{-1}$ when $v=d_1$. Similarly , let $y=s_yA_{-i}$ be the preceding vertex for $q$ when $q\ne d_1$ and let $y=d_{-1}$ when $q=d_1$.
Then $s_v=s_ub$, $s_q=s_yd$ for some $b,d\in T_{-i}$.
We have $hv=s_qcs_v^{-1}s_vA_i=s_qA_i=q$ and $hu=s_ydcb^{-1}s_u^{-1}s_uA_{-i}=s_yA_{-i}=y$.
Thus $h$ takes the edge $(u,v)$ into the edge $(y,q)$ what contradicts the fact that no two distinct edges of $Y$ are $H$-equivalent.
\enddemo

The following statements are obvious corollaries of the properties
of amalgamated free products.
\proclaim {Lemma 1.4}
If $x,y\in G$ and $x=u_1u_2\dots u_k$, $y=v_1\dots v_s$ are their reduced
expressions.
Suppose that $v_s\in A_i$.
Then $xy$ ends in $A_i$ unless $y^{-1}$ is a right segment of $x$.
\endproclaim

\proclaim {Lemma 1.5} If $v$, $u$ are vertices of $Y$, $v\ne d_1$
then $s_u$ is a left segment of $s_v$ if and only if $u\le v$.
\endproclaim

\subhead Transversal elements \endsubhead

\definition{Definition 1.6}
Define the following functions $\rho_i, \sigma_i:\{s_v|v\in VY\}\rightarrow G$, $i=\pm 1$.
If $v=s_vA_i=s_ubA_i\in VY$, where $b\in T_{-i}$ and $u=s_uA_{-i}$ is the preceding vertex of $v$, put $\sigma_{-i}(s_v)=s_u$.
If $v=s_vA_i\in VY_1$, put $\sigma_i(s_v)=s_v$.
If $v=s_vA_i\in V(Y-Y_1)$ and $q=s_qA_i\in VY_1$ is $H$-equivalent to $v$,
then put $\sigma_i(s_v)=s_q$.
Now for $v\in VY_1$ put $\rho_{\pm 1}(s_v)=s_v$.
Suppose $v=s_vA_i\in V(Y-Y_1)$  and $q=s_qA_i\in VY_1$ is $H$-equivalent to $v$. Let $h_v=s_va^{-1}s_q^{-1}$.
Then put $\rho_{-i}(s_v)=s_v$ and $\rho_i(s_v)=s_qa$.

The elements of the set $im(\rho_1)C\cup im(\rho_{-1})C$ are called {\it transversal elements}. 
\enddefinition
This definition is motivated by the work of B.Baumslag [2] who used a similar construction to analyze the subgroup structure of a free product of two groups.
We collect some useful facts about the functions $\rho_i, \sigma_i$ in the following lemma.

\proclaim {Lemma 1.7} Let $v=s_vA_i\in VY$. Then
\roster
\item "(i)" $\sigma_j(s_v)\in Hs_vA_j$, $j=\pm 1$;
\item "(iii)"$\rho_j(s_v)\in Hs_v$, $j=\pm 1$;
\item "(iii)" $\rho_j(s_v)=\sigma_j(s_v)a_j, a_j\in A_j$, $j=\pm 1$;
\item "(iv)" $\sigma_j(s_v)$ is either $1$ or it ends in $A_{-j}$;
\item "(v)" if $v=s_vA_i\in V(Y-Y_1)$ then
$h_v=\rho_{-i}(s_v)\rho_i(s_v)^{-1}$;
\item "(vi)" if $H\sigma_j(s_v)A_j=H\sigma_j(s_w)A_j$ then $\sigma_j(s_v)=\sigma_j(s_w)$, $j=\pm 1$, $v,w\in VY$;
\item "(vii)" if $H\rho_j(s_v)C=H\rho_j(s_w)C$ then $\rho_j(s_v)=\rho_j(s_w)$, $j=\pm 1$, $v,w\in VY$.
\endroster
\endproclaim

\demo {Proof}
Statements (i), (ii), (iii),(iv) and (v) follow immediately from the definitions of $s_v$, $h_v$ $\rho_i$ and $\sigma_i$.

(vi) For any $r\in im(\sigma_i)$ there is $v=s_vA_i\in VY_1$ such that $r=\sigma_i(s_v)=s_v$.

So if $r_1,r_2\in im(\sigma_i)$ and $Hr_1A_i=Hr_2A_i$, let $v_j=s_{v_j}A_i\in VY_1$ be such that $r_j=\sigma_i(s_{v_j})$, $j=1,2$.
Thus $Hs_{v_1}A_i=Hs_{v_2}A_i$. However, no two distinct vertices of $Y_1$ are $H$-equivalent. Therefore $v_1=s_{v_1}A_i=s_{v_2}A_i=v_2$ and $r_1=s_{v_1}=s_{v_2}=r_2$.

(vii) Any element in the image of $\rho_i$ has the form $s_v$, $v\in VY_1$ or $s_va_i$ where $v=s_vA_i\in VY_1$, $a_i\in A_i-C$, $s_va_iA_{-i}\not\in VY_1$.

Suppose $r,y\in im(\rho_i)$ and $HrC=HyC$.
Thus there are $h\in H, c\in C$ such that $h=ycr^{-1}$.
There are several cases to consider.

{\bf Case 1} Assume first $r=s_v,y=s_w$ for some vertices $v,w$ of $Y$.

Suppose first $v$ and $w$ have the same type $A_j$.
Then by Lemma 1.3 $Hs_vC=Hs_wC$ implies $v=w$, $r=s_v=s_w=y$.

Suppose now that $r=s_v$, $v=s_vA_j\in VY$ and $y=s_w$, $w=s_wA_{-j}\in VY$. Since both $r$ and $y$ are in the image of $\rho_i$, one of the vertices $v,w$, say $v$, has type $A_i$ and is in $Y_1$ and the other, say, $w$ has type $A_{-i}$ and is in $Y-Y_1$.
 Consider the preceding vertex $u=s_uA_i\in VY_1$ for $w$. We have $s_w=s_ua_i$ where $a_i\in T_i$.
We know that $h^{-1}=s_vc^{-1}a_i^{-1}s_u^{-1}$.
It is clear that $h^{-1}u=s_vc^{-1}a_i^{-1}s_u^{-1}s_uA_i=s_vA_i=v$.
Since $v,u\in VY_1$, $v=u$ and $s_v=s_u$.
Therefore $h=s_ua_ics_u^{-1}\in H$ and so $a_ic\in A_u$.
This contradicts Lemma 1.2(i).
If $v$ has type $A_{-i}$ then $Hs_vC=Hs_wC$ implies $Hs_vA_{-i}=Hs_wA_{-i}$ that is $Hv=Hw$.
Thus $v$ is the only vertex of $Y_1$ $H$-equivalent to $w$.
Notice that $h(v)=s_wcs_v^{-1}s_vA_{-i}=s_wA_{-i}=w$.
Let $h_w=s_wa^{-1}s_v^{-1}\in H$, $a\in A_{-i}$.
Recall that by Lemma 1.2(ii) $a\not\in C$.
Then $h^{-1}=s_vc^{-1}s_w^{-1}$ and $h^{-1}h_w=s_vc^{-1}a^{-1}s_v^{-1}$.
However $h^{-1}w=v, h_w(v)=w$ and, therefore, $h^{-1}h_w(v)=v$.
Thus $s_vacs_v^{-1}(v)=v$ that is $ac\in A_v$.
But this is impossible by Lemma 1.2(2).

{\bf Case 2} Suppose now that $r=s_{q_1}a_1, y=s_{q_2}a_2$ where $q_1,q_2\in VY_1$ and for some $v_1=s_{v_1}A_i,v_2=s_{v_2}A_i\in V(Y-Y_1)$. $h_{v_1}=s_{v_1}a_1^{-1}s_{q_1}^{-1}$,
$h_{v_2}=s_{v_2}a_2^{-1}s_{q_2}^{-1}$. 
Then $q_1,q_2$ have type $A_i$, $a_1,a_2\in A_i$.
Since $HrC=HyC$, $Hs_{q_1}A_i=Hs_{q_2}A_i$ and therefore $q_1=q_2=q$.
Suppose $r\not=y$. Then  $h_{v_1}hh-{v_2}^{-1}=s_{v_1}cs_{v_2}^{-1}\in H$, and so $Hs_{v_1}C=Hs_{v_2}C$.
Lemma 1.3 implies that $s_{v_1}=s_{v_2}$ and therefore $r=y$.

{\bf Case 3} Suppose now that $r=s_qa$, $y=s_w$, $a\in A_i$, $q=s_qA_i\in VY_1$, $w\in VY$, $h_v=s_va^{-1}s_q^{-1}$, $v\in V(Y-Y_1)$.
Assume first that $w$ has type $A_i$ and therefore $w\in VY_1$.
We have $h=ycr^{-1}=s_wca^{-1}s_q^{-1}\in H$.
Thus $h(q)=s_wca^{-1}s_q^{-1}s_qA_i=s_wA_i=w$.
Therefore $q=w$ since $q,w\in VY_1$.
Hence $h^{-1}=s_qac^{-1}s_q^{-1}\in H$ and $ac^{-1}\in A_q$.
But this contradicts Lemma 1.2(2).
Assume now that $w$ has type $A_{-i}$. Let $u=s_uA_i$ be the preceding vertex of $w$ if $w\ne d_1$ and $u=d_{-1}$ when $w=d_1$. Thus $u\in Y_1$, $s_w=s_ub$, $b\in T_i$.
Then $h=s_ubca^{-1}s_q^{-1}=h\in H$ and so $hq=u$.
This implies $u=q$ since $u,q\in Y_1$.
Therefore $h=s_q(bca^{-1})s_q^{-1}\in H$ and $bca^{-1}\in A_q$.
If $w\ne d_1$ then this contradicts Lemma 1.2(4). If $w=d_1$ then $b=1$ and $ac^{-1}\in A_q$ which contradicts Lemma 1.2(2)
This completes the proof of Lemma 1.7
\enddemo

Lemma 1.7(4) implies that different elements in $im(\rho_j)$ represent different double coset classes $HgC$, $j=\pm 1$. This justifies the term {\it transversal} for the elements of the set $im(\rho_1)C\cup im(\rho_{-1})C$.   
Notice that a left segment of a transversal element is again transversal. 

\definition {Definition 1.8}
Let $g^{-1}$ be a nontransversal element and 
let $w=v_1\dots v_k$, be a reduced form of $g$ with respect to presentation (1).
Let $s\le k$ be the minimal number such that $g=v_1\dots v_s v$ where $v^{-1}$ is a transversal.
We call the expression $v_1\dots v_s$ the {\it nerve} of $w$.
The number $s$ is termed the {\it syllable length of the nerve} of $w$. 
Notice that if $w_1=u_1\dots u_k$ is another reduced form of $g$ and $u_1\dots u_{s_1}$ is the nerve of $w_1$ then $s=s_1$
and $u_1\dots u_sC=v_1\dots v_sC$.

If $g^{-1}$ is transversal and $w=v_1\dots v_k$ is a reduced form of $g$ with respect to presentation (1),
we say that a nerve of $w$ is empty and that it has the syllable length zero.
\enddefinition

\proclaim {Remark} Notice that if $w=v_1\dots v_k$ is a reduced expression with respect to presentation (1), $(v_j\dots v_k)^{-1}$ is a transversal and $(v_{j-1}v_j\dots v_k)^{-1}$ is not a transversal, then $v_1\dots v_{j-1}$ is the nerve of $w$. This immediately follows from the fact that an initial segment of a transversal element is again transversal.
\endproclaim

\proclaim {Lemma 1.9}
\roster
\item "(i)" If $a_v\in A_v-C_v$, $v=s_vA_i\in VY_1$ then $s_va_v$ is not a transversal.
\item "(ii)" If $a_v\in A_v-C_v$,$v=s_vA_i\in VY_1$, $b\in T_i-\{1\}$, $w=s_vbA_{-i}\in VY$, $v$ is a preceding vertex for $w$ and $s_va_vs_v^{-1}$ does not stabilize the edge $(v,w)$  then $s_va_vb$ is not a transversal.
\item "(iii)" If $w=s_wA_i\in V(Y-Y_1)$, $q=s_qA_i\in VY_1$, $a\in A_i$, $b\in A_{-i}-C$ and $h_w=s_wa^{-1}s_q^{-1}$ then $s_qab$ is not a transversal. 
\item "(iv)" Suppose $w=s_wA_i\in V(Y-Y_1)$, $q=s_qA_i\in VY_1$, $a\in A_i$ and $h_w=s_wa^{-1}s_q^{-1}$. Suppose $u=s_uA_{-i}\in VY_1$ is the vertex preceding $w$ and $s_w=s_ub$, $b\in T_{-i}$. Suppose $a_1\in A_i-C$.
Then $s_wa_1$ is not a transversal.
\item "(v)" Suppose $1\ne \rho_i(t)\rho_{-i}(t)=s_va_ia_{-i}^{-1}s_w^{-1}$ where $v=s_vA_i, w=s_wA_{-i}\in VY_1$.
Suppose $b\in A_{i}-C$. Then $s_wa_{-i}b$ is not a transversal.
\item "(vi)" Suppose $q=s_qA_i\in VY_1$, $a_q\in A_q$ and  $a\in A_i$ is such that for some vertex $v\in V(Y-Y_1)$ we have $h_v=s_va^{-1}s_q^{-1}$. Suppose further that $a_qaC\ne aC$.
Then $s_q(a_qa)$ is not a transversal.
\endroster
\endproclaim

\demo {Proof}

(i) Suppose first that $s_va_v\in im(\rho_i)C$.
There are two possibilities.

{\bf Case 1.} There is a vertex $w=s_vbA_{-i}$, $b\in T_i-\{1\}$, such that $v$ precedes $w$ and $a_vC=bC$. Thus $bC\cap A_v\not=\emptyset$ what contradicts Lemma 1.2(1).

{\bf Case 2.} There is a vertex $w=s_wA_i\in Y-Y_1$ $H$-equivalent to $v$, $h_w=s_wa^{-1}s_v^{-1}$ and $a_vC=aC$.
But this contradicts Lemma 1.2(2) which implies $aC\cap A_v=\emptyset$.

Thus $s_va_v$ is not in $im(\rho_i)C$.

Suppose now that $s_va_v\in im(\rho_{-i})C$.
Since $s_va_v$ ends in $A_i$, it means that there is a vertex $w=s_vbA_{-i}$, $b\in T_i$, such that $v$ precedes $w$ and $a_vC=bC$. But this is impossible by Case 1 above.
\smallskip
(ii) Suppose $s_va_vb$ is a transversal.
Assume first that $s_va_vb\in im(\rho_i)C$.
There are two possibilities.

{\bf Case 1.} There is a vertex $u=s_vaA_{-i}$, $a\in T_i$, such that $v$ precedes $v$ and $aC=a_vbC$.
If $a=b$ then $a^{-1}a_va=c\in C$.
Then $h=s_va_vs_v^{-1}\in H$ and $hv=v$.
Moreover, $hw=s_va_vs_v^{-1}s_vaA_{-i}=s_va_vaA_{-i}=s_vacA_{-i}=s_vaA_{-i}=w$. This contradicts our assumption that $h$ does not stabilize the edge $(v,w)$.
If $a\not=b$, $u\not=w$ then by Lemma 1.2(3) $A_vaC\not=A_vbC$.
This contradicts $aC=a_vbC$.

{\bf Case 2.} There is a vertex $y=s_yA_i\in V(Y-Y_1)$ which is $H$-equivalent to $v$ and $h_y=s_ya^{-1}s_v^{-1}$, $a\in A_i$ and $aC=a_vbC$. But Lemma 1.2(4) implies that $A_va\cap bC=\emptyset$ which gives us a contradiction.

Suppose now that $s_va_vb\in im(\rho_{-i})C$.
By Lemma 1.2(1) $bC\cap A_v=\emptyset$, so $a_vb\not\in C$.
Since $s_v(a_vb)$ ends in $A_i$, there is a vertex $u=s_vaA_{-i}\in VY$, $a\in T_i$ such that $v$ precedes $u$ and $aC=a_vbC$. but this is impossible by Case 1 above.
\smallskip
(iii) 
Suppose that $s_qab$ is transversal.
This necessarily implies that $a$ represents the same $C$-coset class as the label of some edge of $Y_1$ emanating from $q$.
But this is impossible by Lemma 1.2(4).
\smallskip
(iv) Suppose $s_wa_1=s_uba_1$ is a transversal. Then $b$ is a label of some edge of $Y_1$ originating from $u$.
This is impossible since the only edge with label $b$ emanating from $u$ is the edge $(u,w)$ and we know that $w\not\in Y_1$.
\smallskip
(v) follows from (iii) and (iv).
\smallskip
(vi) Notice that by Lemma 1.2(2) we have $a_qa\not\in C$.
Assume that $s_q(a_qa)$ is a transversal.
There are two possibilities.

Case 1. There is a positive edge of $Y$ originating from $q$ with label $b\in T_i$ such that $bC=a_qaC$. But by Lemma 1.2(4) we have $bC\cap A_qa=\emptyset$ which gives us a contradiction.

Case 2. There is a vertex $w\in V(Y-Y_1)$ which is $H$-equivalent to $q$ such that $h_w=s_wa_1^{-1}s_q^{-1}$ and $a_qaC=a_1C$. 
Since by assumption $a_qaC\ne aC$, we conclude that $a_1C\ne aC$.
But by Lemma 1.2(3) we have $A_qaC\ne A_qa_1C$ which is impossible. 
\enddemo 

\proclaim {Lemma 1.10}
Suppose $1\ne g=\rho_i(t)\rho_{-i}(t)=s_va_ia_{-i}^{-1}s_q^{-1}$ where $v=s_vA_i, q=s_qA_{-i}\in VY_1$.
Suppose $b\in T_{i}$ is the label of a positive edge of $Y_1$ originating from $v$.
Then $a_{-i}C\ne bC$.
\endproclaim 
\demo{Proof}
There are two cases to consider.

Case 1. Suppose first that $a_i\in T_i$ is the label of an edge $(v,w)\in E^{+}(Y-Y_1)$, $s_w=s_va_i$ and $g=h_w=s_wa_{-i}^{-1}s_q^{-1}$.
Then $a_{-i}C\ne bC$ since $a_i$ and $b$ are the labels of different edges (one is in $Y_1$ and the other is in $Y-Y_1$).

Case 2. Suppose now that $a_{-i}$ is the label of an edge $(q,w)\in E^{+}(Y-Y_1)$, $s_w=s_qa_{-i}$ and $g=h_w^{-1}=(s_wa_{i}^{-1}s_q^{-1})^{-1}$.
Then $bC\cap A_qa_i=\emptyset$ by Lemma 1.2(4) which implies $bC\ne a_iC$.
\enddemo

\subhead Controlling the syllable length of elements of $H$ \endsubhead

Suppose $u=s_vas_v^{-1}\in H$, $a\in A_v-C_v$.
Then $s_vas_v^{-1}$ is a reduced form of $u$ with respect to presentation (1) and we denote it by $w(u)$.
If $1\not=u=\rho_i(t)\rho_{-i}(t)^{-1}=s_va_ia_{-i}^{-1}s_w^{-1}$
then $s_va_ia_{-i}^{-1}s_w^{-1}$ is a reduced form of $u$ and we denote it by $w(u)$.

Recall that each positive edge of $Y$ has a label $a\in T_i$.
Denote $l_{X_i}(a)$ by $l(a)$.
Let $v=s_vA_i\in VY$ and $s_v=a_1\dots a_k$ where $a_j$ is the label of the $j$-th edge of the reduced edge-path from $d_{-i}$ to $v$ in $Y$, $j=1,\dots ,k$.
Then denote $l(a_1)+\dots +l(a_k)$ by $l(s_v)$.
Analogously, for a transversal element $1\ne t=\rho_i(g)\rho_{-i}(g)^{-1}=s_va_ia_{-i}^{-1}s_w^{-1}$ put $l(t)=l(s_v)+l_{X_i}(a_i)+l_{X_{-i}}(a_{-i})+l(s_w)$.
Let $$K=2\underset {t\in im(\rho_{\pm 1})}\to\sum l(t)$$ and $$\Sigma=T_1^{\pm 1}\cup T_{-1}^{\pm 1}\cup \{ a|\  \text{for some}\ v\in V(Y-Y_1)\  \text{ we have}\  h_v=s_va^{-1}s_q^{-1}\}^{\pm 1}.$$

\proclaim {Proposition A}

Let $p$ be a reduced form of $h\in H-C$ with respect to presentation (3). 
Let $U=u_1\dots u_n$ be obtained from $p$ by deleting all those terms which are equal to $1$ in $H$, that is $U$ is a normal form of $h$ with respect to presentation (3). Put $h'=u_1\cdot\dots\cdot u_{n-1}$ when $n>1$.

There is a reduced form with respect to presentation (1) $W=v_1\dots v_m$ of $h$ and, when $n>1$, a reduced form with respect to presentation (1) $W'$ of $h'$ such that the following holds.

\roster 

\item "(i)" If $h\not\in C$ and $u_n=s_vas_v^{-1}$ where $a\in A_v$, $v=s_vA_i$, then $w$ ends in $xs_v^{-1}$ where $x\in A_i$ and  $s_vx^{-1}$ is not a transversal.
If $u_n=h_v^{\pm 1}=\rho_i(g)\rho_{-i}^{-1}(g)=s_va_ia_{-i}^{-1}s_w^{-1}$ then $w$ ends in $xa_{-i}^{-1}s_w^{-1}$ where $x\in A_i-C$ and $s_wa_{-i}x^{-1}$ is not a transversal.
\item "(ii)" The nerve $N$ of $w$ has syllable length not smaller than the nerve $N^{\prime}$ of $w'$.
The syllable length of $N$ is strictly greater than the syllable length of $N^{\prime}$ unless $u_{n-1}=s_vas_v^{-1}$, $a\in A_i$, $v=s_vA_i\in VY_1$ and  $u_n=\rho_i(t)\rho_{-i}(t)^{-1}=s_va_ia_{-i}s_v^{-1}$
(notice that we allow here $s_v=1$).
\item "(iii)" 
For each $u_k=s_vas_v^{-1}$ where $v=s_vA_i\in VY_1$, $a\in A_v$ for some $k=1,\dots ,n$ there is a syllable $v_{i_k}\in A_i-C$ of $w$, called a {\bf core element}, such that 
{\roster
\item if $k_1<k_2$ then $i_{k_1}<i_{k_2}$
\item if $k<n$ then $v_{i_k}=faf'$ where each of $f,f'\in \Sigma$ ;
\item if $k=n$ then $v_{i_k}=fa$ where $f\in \Sigma$;
\item if  $u_k=s_vas_v^{-1}$, $u_{k+1}=s_ubs_u^{-1}=s_ub_1bb_1^{-1}s_u$, where $v$ precedes $u$, $s_u=s_vb_1$, $b_1\in T_i$, $b\in A_{-i}$ then 
$v_{i_k}=fab_1$ and $v_{i_{k+1}}=bf'$ where  $f,f'\in \Sigma$;
\item if  $u_k=s_vas_v^{-1}$, $u_{k+1}=s_ubs_u^{-1}$ where $u$ precedes $v$, $s_v=s_ub_1$, $b_1\in T_{-i}$, $b\in A_{-i}$ then $v_{i_k}=fa$ and $v_{i_{k+1}}=b_1^{-1}bf'$ where $f,f'\in \Sigma$;
\item if $u_k=s_vas_v^{-1}$, $u_{k+1}=s_va_ia_{-i}^{-1}s_w^{-1}$ then $v_{i_k}=faa_i$ where $f\in \Sigma$;
\item if $u_k=s_va_ia_{-i}^{-1}s_w^{-1}$ and $u_{k+1}=s_wbs_w^{-1}$ then $v_{i_{k+1}}=fa_{-i}^{-1}b$ where $f\in \Sigma$;
\item if $u_n=s_vas_v^{-1}$, $v=s_vA_i\in VY_1$, $a\in A_v-C_v$ and $w_n$ ends in $xs_v^{-1}$, $x\in A_i$ then $v_{i_n}=x$;
\item  $v_{k}$ and $v_{k+1}$ are both core elements if and only if
   
either $v_k=v_{i_s}, v_{k+1}=v_{i_{s+1}}$, $u_s=s_vas_v^{-1}$, $u_{s+1}=s_ubs_u^{-1}$ and either $u$ precedes $v$ or $v$ precedes $u$

or $v_k=v_{i_s}$, $v_{k+1}=v_{i_{s+2}}$, $u_s=s_vb_is_v^{-1}$, $u_{s+1}=\rho_i(g)\rho_{-i}^{-1}(g)^{-1}=s_va_ia_{-i}^{-1}s_w^{-1}$, $u_{s+2}=s_wb_{-i}s_w^{-1}$;
\item every $v_k\in A_i$, which is not a core element, has $l_{X_i}(v_k)\le K$.

\endroster}

\endroster 
\endproclaim

\demo {Proof}
We will prove Proposition A by induction on $n$.
Suppose $n=1$.
Recall that $h\not\in C$.

Suppose first $p=(e_1,1,e_2,1,\dots ,1,e_k, u_1 ,e_{k+1},1,\dots e_{k+s},1)$
be a reduced path in ${\Bbb B}$ representing $h$.
Then $s=k$ and $e_{k+i}=e_{k-i+1}^{-1}$, $i=1,\dots ,k$ and $e_i\in Y_1$ for $i=1,\dots k$.Thus $u_1=s_vas_v^{-1}$, $v=s_vA_i\in VY_1$ where $e_1,\dots ,e_k$ is a path in $Y_1$ ending ay $v$ with the label $s_v$.
Notice that $a\not\in C$. Indeed, if $k\ge 1$ and $a=c\in C$ then $u_1=s_vcs_v^{-1}$ stabilizes the edge $e_k$ and $p$ is not a normal form of $h$. If $k=0$ (that is $v=d_1$) and $a=c\in C$ then $h=c\in C$ which, as we assumed, is not the case.
Thus $a\not\in C$ and therefore $s_va^{-1}$ is not a transversal by Lemma 1.9(i).
We have verified statement (i).
Notice that $w=s_vas_v^{-1}$ is the reduced form of $h$ with respect to presentation (1) and the nerve $N$ of $w$ is equal to $s_va$. Put $v_{i_1}=a$ to be the only core element.
Then the rest of the statements of Proposition A are automatically satisfied.

Suppose now that $p=(e_1,\dots ,e_k,e, e_{k+2},\dots ,e_s)$ where
$e_i\in EY_1$, $i=1,\dots ,k,k+2,\dots ,s$ and $e\in E(B-Y_1)$.
Then $u_1=e=\rho_i(t)\rho_{-i}(t)^{-1}=s_va_ia_{-i}^{-1}s_w^{-1}$.
Observe that $s_wa_{-i}a_i^{-1}$ is not a transversal by Lemma 1.9(i).
Observe that $w=s_va_ia_{-i}^{-1}s_w^{-1}$ is the reduced form of $h$ with respect to presentation (1). 
So the nerve $N$ of $w=s_va_ia_{-i}^{-1}s_w^{-1}$ is $N=s_va_i$, there are no core elements and Proposition A for the case $n=1$ is established.
Thus the basis of induction is verified.
\smallskip
Suppose now $n>1$ and Proposition A has been established for smaller values of $n$. There are several cases to consider.
\smallskip
{\bf Case 0.} Suppose that $h'=c\in C$. Therefore $n=2$ and $u_1=c$. Notice that $u_2\not\in A_{\pm 1}$ since $u_1u_2$ is a normal form of $h$ with respect to presentation (3).
\smallskip
{\bf Subcase 0.A.} Suppose that 
$u_2=s_vas_v^{-1}$ where $A_{\pm 1}\ne v=s_vA_i\in VY_1$.

Thus $p=(c,e_1,\dots ,e_k, s_vas_v^{-1}, e_k^{-1},\dots ,e_1^{-1})$ where $k>1$ and $e_1,\dots ,e_k$ is a reduced edge-path in $Y_1$ from $d_1$ to $v$.

Notice that $a\in A_i-C$ since $s_vas_v^{-1}$ does not stabilize the edge $e_k$. Let $s_v=fz$ where $f\in A_j-C$ is the first syllable of $s_v$.
Then $W=(cf)zas_v^{-1}$ is the reduced form of $h$ with respect to presentation (1). The element $s_va^{-1}$ is not a transversal by Lemma 1.9(i). Thus $W$ ends in $s_va^{-1}$ and the nerve $N=(cf)za$ of $W$ has greater syllable length than the nerve $N'=1$ of $W'=c$. Put $v_{i_1}=cf$ and $v_{i_2}=a$ to be the core elements of $w$. Notice that $f\in \Sigma$. All statements of Proposition $A$ are clearly satisfied. 
\smallskip
{\bf Subcase 0.B.} Suppose that $u_2=\rho_i(g)\rho_{-i}(g)^{-1}=s_va_ia_{-i}^{-1}s_w^{-1}\ne 1$ where $v=s_vA_i, w=s_wA_{-i}$, $a_j\in A_j-C$, $j=\pm 1$.

Assume first that $s_v\ne 1$. Then $s_v=fz$ where $f\in A_j-C$ is the first syllable of $s_v$. Then
$W=(cf)za_ia_{-i}^{-1}s_w^{-1}$ is a reduced form of $h$ with respect to presentation (1). By Lemma 1.9(v) the element $s_wa_{-i}a_i^{-1}$ is not a transversal. Therefore the nerve $N$ of $W$ is equal to $N=(cf)za_i$ and it has reater syllable length than the nerve $N'=1$ of $W'=c$. Put $v_{i_1}=cf$ to be the only core element for $W$. Notice that $f\in \Sigma$.
All statements of Proposition $A$ are clearly satisfied. 

Suppose now that $s_v=1$. Then $W=(ca_i)a_{-i}^{-1}s_w^{-1}$ is the reduced form of $h$ with respect to presentation (1).  By Lemma 1.9(v) the element $s_wa_{-i}a_i^{-1}c^{-1}$ is not a transversal. Therefore the nerve $N$ of $W$ is equal to $N=(ca_i)$ and it has reater syllable length than the nerve $N'=1$ of $W'=c$. Put $v_{i_1}=ca_i$ to be the only core element for $W$. Notice that $a_i\in \Sigma$. All statements of Proposition $A$ are clearly satisfied.

\smallskip
{\bf Case 1.} Suppose that $h'\not\in C$, $u_{n-1}=s_vb_0s_v^{-1}$, $v=s_vA_i$, $b_0\in A_v$, 
$u_n=s_was_w^{-1}$, $a\in A_v$, $w=s_wA_j\in VY_1$ and $w\ne=v$.

Thus $p=(g_1,e_1,\dots ,g_k=s_vb_0s_v^{-1}, e_k,1,e_{k+1},1,\dots 1,e_l, g_{l+1}=s_was_w^{-1},e_{l+1},1,\dots ,e_r,1)$ where
$e_i\in EY_1$ for $i=k,\dots ,r$, $e_k,e_{k+1},\dots ,e_r$ is a path in $Y_1$ from $v$ to $d_1$, $e_{l+1},\dots ,e_r$ is a reduced path in $Y_1$ from $w$ to $d_1$. Let $e_k,\dots ,e_l={\hat z}^{-1} {\hat u}$, $e_{l+1},\dots ,e_r={\hat u}^{-1}{\hat y}^{-1}$ where
${\hat u}^{-1}$ is the maximal initial segment of $e_{l+1},\dots ,e_r$ which is
cancelled in $e_k,\dots ,e_l,e_{l+1},\dots ,e_r$.
Thus ${\hat z}^{-1}{\hat y}^{-1}$ is a reduced path in $Y_1$ from $v$ to $d_1$.
Let ${\hat z}^{-1}{\hat y}^{-1}=e_{k'},\dots ,e_{t'}$. Then $p'=(g_1,e_1,\dots ,g_k=s_vb_0s_v^{-1},e_{k'},1,\dots ,1,e_{t'},1)$ is a reduced form of $h^{\prime}=u_1\dots u_{n-1}$ with respect to presentation (3). Therefore by induction $W'=pxs_v^{-1}$ where $x\in A_i$, $s_vx^{-1}$ is not a transversal and $N^{\prime}=px$ is the nerve of $W'$. Also by induction we know that $v_{i_{n-1}}^{\prime}=x$ is the last core element of $p'$ and that $x=\sigma b_0$ for some $\sigma\in \Sigma$. Let $u$ be the label of $\hat u$, $y$ be the label of $\hat y$ and $z$ be the label of $\hat z$.
Therefore $s_v=yz$ and $s_w=yu$.

\smallskip
{\bf Subcase 1.A.} Suppose first that both $\hat z$ and $\hat u$ are nonempty.

  Then $s_v=yq$ and $s_w=yu$ and $u_n=s_was_w^{-1}=yuau^{-1}y^{-1}$.
Notice that $a\not\in C$ since if $a\in C$ then $u_n$ fixes the last edge of $\hat u$ which contradicts our assumption that $p$ is a reduced form for $h$ with respect to presentation (3).
Thus $h=pxz^{-1}y^{-1}yuas_w^{-1}=pxz^{-1}uas_w^{-1}$.
Suppose that $\hat y$ ends in a vertex of type $A_k$.
Let $z=f_1z_1$ where $f_1\in T_k$ be the label of the first edge of $\hat z$. Let $u=f_2u_1$ where where $f_2\in T_k$ be the label of the first edge of $\hat u$.
Clearly $f_1C\ne f_2C$ by definition of $\hat u$ and $\hat z$ and so $f_1^{-1}f_2\not\in C$.
Therefore
$W=pxz_1^{-1}(f_1^{-1}f_2)u_1as_w^{-1}$ is a normal form for $h$ with respect to presentation (1).
It is clear that $s_w$ is a transversal.
Besides $s_wa^{-1}$ is not a transversal by Lemma 1.9(i).
Thus the nerve $N$ of $W$ is equal to $pxz_1^{-1}(f_1^{-1}f_2)u_1a$ and it has greater syllable length
then the nerve $N^{\prime}=px$ of $W'$.
Now take the set of core element of $W'$, add to it $v_{i_n}=a$ and declare the result to be the set of core elements of $w_n$.
All statements of Proposition A are clearly satisfied by induction.

\smallskip
{\bf Subcase 1.B.} Suppose that $\hat z$ is empty and $\hat u$ is nonempty.

Then $s_v=y$ and $s_w=yu$ and so $h=pxs_v^{-1}\cdot s_was_w^{-1}=pxy^{-1}yuas_w^{-1}=pxuas_w^{-1}$.
Let $u=fu_1$ where $f$ is the label of the first edge of $u$.
Thus $x,f\in A_i$ since $v$ is the vertex of type $A_i$.
Notice that $x\cdot f\not\in C$.
Indeed, if $xf=c\in C$ then $s_vx^{-1}=s_vfc^{-1}$ which is a transversal element.
This clearly contradicts the inductive assumption that $s_vx^{-1}$ is not a transversal. Observe also that $a\not\in C$ since if $a\in C$ then $u_n$ fixes the last edge of $\hat u$ which contradicts our assumption that $p$ is a reduced form for $h$ with respect to presentation (3).
Thus $W=p(xf)u_1as_w^{-1}$  is a normal form for $h$ with respect to presentation (1).
Again we see that $s_w$ is a transversal and $s_wa^{-1}=yua^{-1}$ is not a transversal by Lemma 1.9(i).
Thus the nerve $N$ of $W$ is equal to $p(xf)u_1a$ and it has greater syllable length
then the nerve $N^{\prime}=px$ of $W'$.
Recall that the last core element of $W'$ is $v_{i_{n-1}}^{\prime}=x=\sigma b_0$where $\sigma\in \Sigma$. 
Now take the set of core element of $W'$, replace $v_{i_{n-1}}^{\prime}=\sigma b_0$ by $v_{i_{n-1}}=\sigma b_0f$, add $v_{i_n}=a$ and declare the result to be the set of core elements of $W$.
All statements of Proposition A are clearly satisfied by induction.
\smallskip
{\bf Subcase 1.C.} Suppose that $\hat z$ is nonempty and $\hat u$ is empty.

Then $s_v=yz$ and $s_w=y$.
In this case $h=pxs_v^{-1}s_was_w^{-1}=pxz^{-1}y^{-1}yay^{-1}=pxz^{-1}ay^{-1}$.
Since $p$ is a reduced form for $h$ with respect to presentation (3), the element $s_was_w^{-1}$ does not stabilize the first edge of $z$.
Thus if $f\in T_j$ is the label of this edge and $z=fz_1$, then
$f^{-1}a\in A_j-C$ by Lemma 1.2(1) and $s_w^{-1}a^{-1}f=y^{-1}a^{-1}f$ is not a transversal by Lemma 1.9(ii).
Therefore 
$W=pxz_1^{-1}(f^{-1}a)y^{-1}$ is a normal form for $h$ with respect to presentation (1).
Since $y=s_w$ is a transversal and $s_w^{-1}a^{-1}f=y^{-1}a^{-1}f$ is not a transversal, we conclude that the nerve $N$ of $W$
is $pxz_1^{-1}(f^{-1}a)$ and it has greater syllable length than $N^{\prime}=px$.
Now take the set of core element of $w'$, add to it $v_{i_n}=f^{-1}a$ and declare the result to be the set of core elements of $w_n$. 
All statements of Proposition A are clearly satisfied by induction.
\smallskip
{\bf Subcase 1.D}  Suppose that both $\hat z$ and $\hat u$ are empty. Then $w=v$ which contradicts our assumptions.
\smallskip

{\bf Case 2.} Suppose that $h'\not\in C$, $u_{n-1}=s_vb_0s_v^{-1}$, $v=s_vA_i$, $b_0\in A_v$, 
$1\ne u_n=\rho_j(g)\rho_{-j}(g)^{-1}=s_tb_jb_{-j}^{-1}s_q^{-1}$ where $q=s_qA_{-j}, t=s_tA_{j} \in VY_1$ and $b_{\pm j}\in A_{\pm j}-C$.

Thus $p=(g_1,e_1,\dots ,g_k=s_vb_0s_v^{-1}, e_k,1,e_{k+1},1,\dots 1,e_l, g_{l+1}=s_tb_jb_{-j}^{-1}s_q^{-1}, e_{l+1},1,\dots ,e_r,1)$ where
$e_i\in EY_1$ for $i=k,\dots ,r$, $e_k,e_{k+1},\dots ,e_l$ is a path in $Y_1$ from $v$ to $t$, $e_{l+1},\dots ,e_r$ is a reduced path in $Y_1$ from $q$ to $d_1$. Let $d_1,\dots d_s$ be the reduced path in $Y_1$ from $t$ to $d_1$. 

Then $e_k,e_{k+1},\dots ,e_l,d_1,\dots ,d_s$ is a path in $Y_1$ from $v$ to $d_1$.
Let $e_k,e_{k+1},\dots ,e_l={\hat z}^{-1} {\hat u}$, $d_1,\dots ,d_s={\hat u}^{-1}{\hat y}^{-1}$ where
${\hat u}^{-1}$ is the maximal initial segment of $d_1,\dots ,d_s$ which is cancelled in the product $e_k,e_{k+1},\dots ,e_l,d_1,\dots ,d_s$. 
Then ${\hat u}^{-1}{\hat y}^{-1}=e_k'\dots e_m'$ is a reduced path in $Y_1$ from $v$ to $d_1$.  
Therefore $p'=(g_1,e_1,\dots ,g_k=s_vb_0s_v^{-1},e_{k}',1,\dots ,1,e_{m}',1)$ is a reduced form of $h^{\prime}=u_1\dots u_{n-1}$ with respect to presentation (3).

By induction $W'=pxs_v^{-1}$, where $x\in A_i$, $s_vx^{-1}$ is not a transversal and $N^{\prime}=px$ is the nerve of $W'$. Also by induction we know that for some $\sigma\in \Sigma$ $v_{i_{n-1}}^{\prime}=x=\sigma b_0$ is the last core element of $W'$.
Denote the labels of $\hat u, \hat z, \hat y$ by $u, z, y$.
Therefore $s_v=yz$ and $s_t=yu$.

\smallskip
{\bf Subcase 2.A.} Suppose that $\hat u$ is empty and $\hat z$ is non-empty.

Then $s_v=yz$ and $s_t=y$.
Therefore $h=pxs_v^{-1}s_tb_jb_{-j}^{-1}s_q^{-1}=pxz^{-1}y^{-1}yb_jb_{-j}^{-1}s_q^{-1}=pxz^{-1}b_jb_{-j}^{-1}s_q^{-1}$.
Let $z=f_1z_1$ where $f_1$ is the label of the first edge of $z$.
Then $f_1,b_j\in A_{j}$ and either $b_j$ is a label of the edge $(t,w)\in V(Y-Y_1)$ and $h=h_w$ or $b_{-j}$ is a label of the edge $(q,w)\in V(Y-Y_1)$ and $h=h_w^{-1}$.
In the first case $f_1C\ne fC$ since the first edge of $z$ is in $Y_1$ and $(t,w)\in E(Y-Y_1)$. In the second case $f_1C\cap A_tb_j=\emptyset$ by Lemma 1.2(4).
Thus $(f_1^{-1}b_j)\not\in C$ and
$W=pxz_1^{-1}(f_1^{-1}b_j)b_{-j}^{-1}s_q^{-1}$ is the normal form of $h$ with respect to presentation (1).
Notice that $s_qb_{-j}$ is transversal and $s_qb_{-j}(b_j^{-1}f_1)$ is not transversal by Lemma 1.9(v).
Thus the nerve $N$ of $W$ is $pxz_1^{-1}(f_1^{-1}b_j)$ and it has greater syllable length than the nerve $N'=px$ of $W'$.
Take the core elements of $W'$ and declare them to be the core elements of $W$.
Proposition A now follows from the inductive hypothesis.
\smallskip
{\bf Subcase 2.B} Suppose that $\hat u, \hat z$ are empty.

Then $v=t$, $i=j$, $s_v=y$ and $s_t=y$.
Therefore $h=pxs_v^{-1}s_tb_jb_{-j}^{-1}s_q^{-1}=pxy^{-1}yb_jb_{-j}^{-1}s_q^{-1}=pxb_jb_{-j}^{-1}s_q^{-1}$.
Note that $x,b_j\in A_{j}$, $b_{-j}\in A_{-j}$.
Observe that $xb_j\not\in C$ since if $xb_j=c\in C$ then $s_vx^{-1}=s_tx^{-1}=s_tb_jc^{-1}$ is a transversal which contradicts our assumptions. Recall also that $b_{-j}\not\in C$.
Thus $W=p(xb_j)b_{-j}s_q^{-1}$ is is the normal form of $h$ with respect to presentation (1).
Lemma 1.9(v) implies that $g=s_qb_{-j}(b_j^{-1}x^{-1})$ is not a transversal.
Thus the nerve $N$ of $W$ is equal to $p(xb_j)$ and it has the same syllable length as the nerve $N'=px$ of $W'$.
By the inductive hypothesis $v_{i_{n-1}}^{\prime}=x=\sigma b_0$ is the last core element of $W'$ for some $\sigma\in \Sigma$.
We take the collection of core elements of $W'$ replace $v_{i_{n-1}}^{\prime}=\sigma b_0$ by $v_{i_{n-1}}=\sigma b_0b_j$ and declare this to be the collection of core elements of $W$.
All statements of Proposition A are clearly satisfied by induction.
\smallskip
{\bf Subcase 2.C}
Suppose now that $\hat u$ and $\hat z$ are nonempty.

Then $s_v=yz$ and $s_t=yu$. Therefore
$h=pxs_v^{-1}s_tb_jb_{-j}^{-1}s_q^{-1}=pxz^{-1}y^{-1}yub_jb_{-j}^{-1}s_q^{-1}=pxz^{-1}ub_jb_{-j}^{-1}s_q^{-1}$.
Assume that $\hat y$ is a path from $d_1$ to the vertex of type $A_k$.
Let $z=f_1z_1$ where $f_1\in T_k$ is the label of the first edge of $\hat z$ and let $u=f_2u_1$ where $f_2\in T_k$ is the label of the first edge of $\hat u$.
Notice that $f_1C\ne f_2C$ by definition of $\hat u$ and $\hat z$.
Thus $(f_1^{-1}f_2)\not\in C$ and
$W=pxz_1^{-1}(f_1^{-1}f_2)u_1b_jb_{-j}^{-1}s_q^{-1}$ is the normal form for $h$ with respect to presentation (1).
Again we observe that $s_qb_{-j}$ is a transversal and $s_qb_{-j}b_j^{-1}$ is not a transversal by Lemma 1.9(v).
Thus the nerve of $W$ is equal to $N=pxz_1^{-1}(f_1^{-1}f_2)u_1b_j$ and it has greater syllable length than the nerve $N'=px$ of $W'$.
Take the core elements of $W'$ and declare them to be the core elements of $W$.
It is clear that all statements of Proposition A follow from the inductive hypothesis.
\smallskip
{\bf Subcase 2.D} Suppose that $\hat u$ is nonempty and $\hat z$ is empty.

Then $s_v=y$, $v=yA_i$ and $s_t=yu$.
Therefore $h=pxs_v^{-1}s_tb_jb_{-j}^{-1}s_q^{-1}=pxy^{-1}yub_jb_{-j}^{-1}s_q^{-1}=pxub_jb_{-j}^{-1}s_q^{-1}$.
Let $u=f_1u_1$ where $f_1\in T_i$ is the label of the first edge of $u$.
Then $(xf_1)\not\in C$. Indeed, if $xf_1=c\in C$ then $s_vx^{-1}=s_vf_1c^{-1}$ is a transversal which contradicts our assumptions.
Thus $W=p(xf_1)u_1b_jb_{-j}^{-1}s_q^{-1}$ is the normal form for $h$ with respect to presentation (1).
As in the previous case $s_qb_{-j}$ is a transversal and $s_qb_{-j}b_j^{-1}$ is not a transversal by Lemma 1.9(v).
So the nerve $N$ of $W$ is $p(xf_1)u_1b_j$ and it has greater syllable length than the nerve $N'=px$ of $W'$.

Recall that by inductive hypothesis $v_{i_{n-1}}^{\prime}=x=\sigma b_0$ is the last core element of $W'$.
Take the core elements of $W'$ and replace $v_{i_{n-1}}^{\prime}=\sigma b_0$ by $v_{i_{n-1}}=\sigma b_0f_1$ to get the collection of core elements of $W$. Proposition A follows now from the inductive hypothesis.

\smallskip

{\bf Case 3}. $h'\not\in C$, $1\ne u_{n-1}=\rho_i(g)\rho_{-i}(g)^{-1}=s_va_ia_{-i}^{-1}s_w^{-1}$ and $1\ne u_{n}=\rho_j(g')\rho_{-j}(g')^{-1}=s_tb_jb_{-j}^{-1}s_q^{-1}$  where $v=s_vA_i, w=s_wA_{-i}, t=s_tA_j, q=s_q^{-j}\in VY_1$, $a_\pm i\in A_{\pm i}-C$, $b_{\pm j}\in A_{\pm j}-C$.

Thus $p=(g_1,e_1,\dots ,g_k=s_va_ia_{-i}^{-1}s_w^{-1},e_k,1,\dots ,1,e_l,g_{l+1}=s_tb_jb_{-j}^{-1}s_q^{-1}, e_{l+1},\dots ,e_r,1)$ where $e_i\in EY_1$ for $i\ge k$, $e_k,\dots ,e_l$ is a path in $Y_1$ from $w$ to $t$ and $e_{l+1},\dots ,e_r$ is a reduced path in $Y_1$ from $q$ to $d_1$.
Let $d_1,\dots ,d_s$ be the reduced path in $Y_1$ from $t$ to $d_1$.  Then $e_k,\dots ,e_l,d_1,\dots ,d_s$ is a path in $Y_1$ from $w$ to $d_1$.
Let $e_k,\dots ,e_l={\hat z}^{-1}u$ and $d_1,\dots ,d_s={\hat u}^{-1}{\hat y}^{-1}$ where $\hat u$ is the maximal terminal segment of $e_k,\dots ,e_l$ which is cancelled in $e_k,\dots ,e_l,d_1,\dots ,d_s$.
Then ${\hat u}^{-1}{\hat y}^{-1}=e_k'\dots e_m'$ is a reduced path in $Y_1$ from $w$ to $d_1$.  
Therefore $p'=(g_1,e_1,\dots ,g_k=s_va_ia_{-i}^{-1}s_w^{-1},e_{k}',1,\dots ,1,e_{m}',1)$ is a reduced form of $h^{\prime}=u_1\dots u_{n-1}$ with respect to presentation (3).
By induction $W'=pxa_{-i}^{-1}s_w^{-1}$, where $x\in A_i$, $s_wa_{-i}x^{-1}$ is not a transversal and $N^{\prime}=px$ is the nerve of $W'$. 
Denote the labels of $\hat u, \hat z, \hat y$ by $u, z, y$.
Thus $s_w=yz$ and $s_t=yu$.

\smallskip
{\bf Subcase 3.A.} Suppose that $\hat z$ is nonempty and $\hat u$ is empty.

Then $s_w=yz$, $s_t=y$. Therefore $h=pxa_{-i}^{-1}s_w^{-1}s_tb_jb_{-j}^{-1}s_q^{-1}=pxa_{-i}^{-1}z^{-1}y^{-1}yb_jb_{-j}^{-1}s_q^{-1}=pxa_{-i}^{-1}z^{-1}b_jb_{-j}^{-1}s_q^{-1}$.
Let $z=f_1z_1$ where $f_1\in T_j$ is the label of the first edge of $z$.
Observe that $f_1^{-1}b_j\not\in C$ by Lemma 1.10 and therefore
$W=pxa_{-i}^{-1}z_1^{-1}(f_1^{-1}b_j)b_{-j}^{-1}s_q^{-1}$ is the normal form of $h$ with respect to presentation (1).
The element $s_qb_{-j}$ is transversal and $s_qb_{-j}(f_1^{-1}b_j)^{-1}$ is not transversal by Lemma 1.9(v). Therefore the nerve $N$ of $W$ is equal to $pxa_{-i}^{-1}z_1^{-1}(f_1^{-1}b_j)$ and it has greater syllable length than the nerve $N'=px$ of $W'$.
We take the collection of core elements of $W'$ and declare them to be the core elements of $W$.
Proposition A follows now from the inductive hypothesis.

{\bf Subcase 3.B.} Suppose that $\hat u$ is nonempty and $\hat z$ is empty.

Then $s_w=y$ and $s_t=yu$. We have
$h=pxa_{-i}^{-1}s_w^{-1}s_tb_jb_{-j}^{-1}s_q^{-1}=pxa_{-i}^{-1}y^{-1}yub_jb_{-j}^{-1}s_q^{-1}=pxa_{-i}^{-1}ub_jb_{-j}^{-1}s_q^{-1}$.
Let $u=f_1u_1$ where $f_1\in T_{-i}$ is the label of the first edge of $\hat u$. Then $(a_{-i}^{-1}f_1)\not\in C$ by Lemma 1.10. Thus $W=px(a_{-i}^{-1}f_1)u_1b_jb_{-j}^{-1}s_q^{-1}$ is the normal form of $h$ with respect to presentation (1).
The element $s_qb_{-j}$ is transversal and $s_qb_{-j}b_j^{-1}$ is not transversal by Lemma 1.9(v). Therefore the nerve $N$ of $W$ is equal to $px(a_{-i}^{-1}f_1)u_1b_j$ and it has greater syllable length than the nerve $N'=px$ of $W'$.
We take the collection of core elements of $W'$ and declare them to be the core elements of $W$.
Proposition A follows now from the inductive hypothesis.
\smallskip
{\bf Subcase 3.C.} Suppose that $\hat u$ and $\hat z$ are nonempty.

Then $s_w=yz$ and $s_t=yu$ and $h=pxa_{-i}^{-1}s_w^{-1}s_tb_jb_{-j}^{-1}s_q^{-1}=pxa_{-i}^{-1}z^{-1}y^{-1}yub_jb_{-j}^{-1}s_q^{-1}=pxa_{-i}^{-1}z^{-1}ub_jb_{-j}^{-1}s_q^{-1}$.
Assume that $\hat y$ ends in a vertex of type $A_k$.
Let $z=f_1z_1$ and $u=f_2u_1$ where $f_1\in T_k$ is the label of the first edge of $\hat z$ and $f_2\in T_k$ is the label of the first edge of $\hat u$. Clearly $f_1C\ne f_2C$ and so $f_1^{-1}f_2\not\in C$.
Thus $W=pxa_{-i}^{-1}z_1^{-1}(f_1^{-1}f_2)u_1b_jb_{-j}^{-1}s_q^{-1}$ is the normal form of $h$ with respect to presentation (1).
The element $s_qb_{-j}$ is transversal and $s_qb_{-j}b_j^{-1}$ is not transversal by Lemma 1.9(v).
Therefore the nerve $N$ of $W$ is equal to $pxa_{-i}^{-1}z_1^{-1}(f_1^{-1}f_2)u_1b_j$ and it has greater syllable length than the nerve $N'=px$ of $W'$.
We take the collection of core elements of $W'$ and declare them to be the core elements of $W$.
Proposition A follows now from the inductive hypothesis.

\smallskip
{\bf Subcase 3.D.} Suppose that $\hat u$ and $\hat z$ are empty.

Then $s_t=s_w=y$, $t=w$ and $-i=j$.
We have $h=pxa_{-i}^{-1}s_w^{-1}s_tb_jb_{-j}^{-1}s_q^{-1}=pxa_{-i}^{-1}y^{-1}yb_jb_{-j}^{-1}s_q^{-1}=pxa_{-i}^{-1}b_jb_{-j}^{-1}s_q^{-1}$.
By Lemma 1.7(vii) either $a_{-i}^{-1}b_j\not\in C$ or $u_n=u_{n-1}^{-1}$. The later is impossible since $p$ is the reduced form for $h$ with respect to presentation (3).
Thus  $a_{-i}^{-1}b_j\not\in C$ and 
$W=px(a_{-i}^{-1}b_j)b_{-j}^{-1}s_q^{-1}$ is the normal form of $h$ with respect to presentation (1).
The element $s_qb_{-j}$ is transversal and $s_qb_{-j}(a_{-i}^{-1}b_j)^{-1}$ is not transversal by Lemma 1.9(v).
Therefore the nerve $N$ of $W$ is equal to $px(a_{-i}^{-1}b_j)$ and it has greater syllable length than the nerve $N'=px$ of $W'$.
We take the collection of core elements of $W'$ and declare them to be the core elements of $W$.
Proposition A follows now from the inductive hypothesis.

\smallskip
{\bf Case 4.} Suppose that $h'\not\in C$,  $1\ne u_{n-1}=\rho_i(g)\rho_{-i}(g)^{-1}=s_va_ia_{-i}^{-1}s_w^{-1}$ and $u_{n}=s_tbs_t^{-1}$  where $v=s_vA_i, w=s_wA_{-i}, t=s_tA_j\in VY_1$, $a\pm i\in A_{\pm i}-C$, $b\in A_{j}$.

Thus $p=(g_1,e_1,\dots ,g_k=s_va_ia_{-i}^{-1}s_w^{-1},e_k,1,\dots ,1,e_l,g_{l+1}=s_tbs_t^{-1}, e_{l+1},\dots ,e_r,1)$ where $e_i\in EY_1$ for $i\ge k$, $e_k,\dots ,e_l$ is a path in $Y_1$ from $w$ to $t$ and $e_{l+1},\dots ,e_r$ is a reduced path in $Y_1$ from $t$ to $d_1$.
Then $e_k,\dots ,e_l,e_{l+1},\dots ,e_r$ is a path in $Y_1$ from $w$ to $d_1$.
Let $e_k,\dots ,e_l={\hat z}^{-1}u$ and $d_1,\dots ,d_s={\hat u}^{-1}{\hat y}^{-1}$ where $\hat u$ is the maximal terminal segment of $e_k,\dots ,e_l$ which is cancelled in $e_k,\dots ,e_l,e_{l+1},\dots ,e_r$.
Then ${\hat u}^{-1}{\hat y}^{-1}=e_k',\dots ,e_m'$ is a reduced path in $Y_1$ from $w$ to $d_1$.  
Therefore $p'=(g_1,e_1,\dots ,g_k=s_va_ia_{-i}^{-1}s_w^{-1},e_{k}',1,\dots ,1,e_{m}',1)$ is a reduced form of $h^{\prime}=u_1\dots u_{n-1}$ with respect to presentation (3).
By induction $W'=pxa_{-i}^{-1}s_w^{-1}$, where $x\in A_i$, $s_wa_{-i}x^{-1}$ is not a transversal and $N^{\prime}=px$ is the nerve of $W'$. 
Denote the labels of $\hat u, \hat z, \hat y$ by $u, z, y$.
Thus $s_w=yz$ and $s_t=yu$.
\smallskip
{\bf Subcase 4.A.} Suppose that $\hat z$ is nonempty and $\hat u$ is empty.

Then $s_w=yz$, $s_t=y$, $t=yA_j$. Therefore 
$h=pxa_{-i}^{-1}s_w^{-1}s_tbs_t^{-1}=pxa_{-i}^{-1}z^{-1}y^{-1}ybs_t^{-1}=pxa_{-i}^{-1}z^{-1}bs_t^{-1}$.
Notice that $\hat z$ starts at $t=yA_j$.
Let $z=f_1z_1$, where where $f_1\in T_j$ is the label of the first edge of $\hat z$.
Then $f_1^{-1}b\not\in C$ by Lemma 1.2(1).
Thus $W=pxa_{-i}^{-1}z_1^{-1}(f_1^{-1}b)s_t^{-1}$ is the normal form of $h$ with respect to presentation (1).
The element $s_q$ is transversal and $s_tbs_t^{-1}$ does not stabilize the first edge of $\hat z$ since $p$ is a reduced form of $h$ with respect to presentation (3).
Therefore $s_qb^{-1}f_1$ is not transversal by Lemma 1.9(ii).
Thus the nerve $N$ of $W$ is $pxa_{-i}^{-1}z_1^{-1}(f_1^{-1}b)$ and it has greater syllable length than the nerve $N'=px$ of $W'$.
Take the core elements of $W'$, add to them $v_{i_n}=f_1^{-1}b$ and declare the result the collection of core elements of $W$.
Proposition A follows now from the inductive hypothesis.

\smallskip
{\bf Subcase 4.B.} Suppose that $\hat u$ is nonempty and $\hat z$ is empty.

Then $s_w=y$, $w=yA_{-i}$, $s_t=yu$. Therefore 
$h=pxa_{-i}^{-1}s_w^{-1}s_tbs_t^{-1}=pxa_{-i}^{-1}y^{-1}yubs_t^{-1}=pxa_{-i}^{-1}ubs_t^{-1}$.
Notice that $\hat u$ starts at $w$ and ends at $t$.
Let $u=f_1u_1$, where where $f_1\in T_{-i}$ is the label of the first edge of $\hat z$.
Then $a_{-i}^{-1}f_1\not\in C$ by Lemma 1.9(v).
Thus $W=px(a_{-i}^{-1}f_1)u_1bs_t^{-1}$ is the normal form of $h$ with respect to presentation (1).
Since $\hat u$ is nonempty, the element $s_tbs_t^{-1}$ does not fix the last edge of $\hat u$ because $p$ is the reduced form of $h$. Therefore $b\in A_t-C_v$.
This implies that $s_vb^{-1}$ is not a transversal by Lemma 1.9(i).
Thus the nerve $N$ of $W$ is $px(a_{-i}^{-1}f_1)u_1b$ and it has greater syllable length than the nerve $N'=px$ of $W'$.
Take the core elements of $W'$, add to them $v_{i_n}=b$ and declare the result the collection of core elements for $W$.
Proposition A follows now from the inductive hypothesis.

\smallskip
{\bf Subcase 4.C.} Suppose that $\hat z$, $\hat u$ are empty.

Then $s_w=y=s_t$, $-i=j$, $w=yA_{-i}=yA_j=t$.
Therefore 
$h=pxa_{-i}^{-1}s_w^{-1}s_tbs_t^{-1}=pxa_{-i}^{-1}y^{-1}ybs_t^{-1}=pxa_{-i}^{-1}bs_t^{-1}$. 

Suppose $a_{-i}^{-1}b\in C$. There are two possibilities.
First, it can happen that $a_{-i}$ is the label of an edge originating from $t$. This is clearly impossible since Lemma 1.2(1) implies $a_{-i}C\cap A_t=\emptyset$.
Secondly, it is possible that $u_{n-1}=h_{w'}=s_{w'}a_{-i}^{-1}s_t^{-1}$ where $w'\in V(Y-Y_1)$ is some vertex $H$-equivalent to $t=w$. Recall that $b\in A_t=A_w$.
Then $A_ta_{-i}\cap C=\emptyset$ by Lemma 1.2(2) and so $a_{-i}^{-1}b\not\in C$ which gives us a contradiction.
Thus $W=px(a_{-i}^{-1}b)s_t^{-1}$ is the normal form of $h$ with respect to presentation (1).

Suppose now that $s_t(b^{-1}a_{-i})$ is a transversal.
There are again two possibilities to consider.

First, suppose that $u_{n-1}=h_{w'}^{-1}=s_va_is_{w'}^{-1}=s_va_ia_is_w^{-1}$ for a vertex $w'\in V(Y-Y_1)$. Then $a_{-i}$ is the label of the edge $(w,w')=(t,w')\in E(Y-Y_1)$.
Since $s_t(b^{-1}a_{-i})$ is a transversal, Lemma 1.9(ii) implies that $s_tb^{-1}s_t^{-1}$ (and so $s_tbs_t^{-1}$) stabilizes the edge $(t,w')=(w,w')$. 
Recall that $h_{w'}$ conjugates the subgroup $s_{w'}(A_{w'}\cap C)s_{w'}^{-1}$ into a subgroup of $s_vA_vs_v^{-1}$. 
Therefore $b=a_{-i}ca_{-i}^{-1}$ and
$h_{w'}^{-1}s_tbs_t^{-1}h_{w'}=s_va_vs_v^{-1}$ for some $a_v\in A_v$.
Thus $u_{n-1}u_n=s_va_vs_v^{-1}u_n$ which contradicts the fact that $p$ is a reduced form for $h$ with respect to presentation (3).
Therefore in this case $s_t(b^{-1}a_{-i})$ is not a transversal. 

Secondly, suppose that $u_{n-1}=h_{w'}=s_{w'}a_{-i}^{-1}s_t^{-1}=s_va_ia_{-i}^{-1}s_t^{-1}$ where $w'\in V(Y-Y_1)$ is some vertex $H$-equivalent to $t=w$. 
Then $s_{w'}=s_va_i$ and $a_i$ is the label of the edge $(v,w')\in E(Y-Y_1)$.
Recall that $b\in A_t=A_w$. Since $s_t(b^{-1}a_{-i})$ is a transversal, Lemma 1.9(vi) implies that $b^{-1}a_{-i}=a_{-i}c$ for some $c\in C$.
Thus $b\in A_t\cap a_{-i}Ca_{-i}^{-1}$. Recall that in this situation $h_{w'}^{-1}$ conjugates the subgroup $A_t\cap a_{-i}Ca_{-i}^{-1}$ of $A_t$ into the subgroup $s_{w'}(A_{w'}\cap C)s_{w'}^{-1}$.
Thus $h_{w'}(s_tbs_t^{-1})h_{w'}^{-1}=s_v(a_ic_1a_i^{-1})s_v^{-1}=s_va_vs_v^{-1}$. Consequently, we have $u_{n-1}u_n=s_va_vs_v^{-1}\cdot u_{n-1}$ which contradicts the fact that $p$ is a reduced form for $h$ with respect to presentation (3).
Therefore in this case $s_t(b^{-1}a_{-i})$ is not a transversal. 

We have established that that $s_t(b^{-1}a_{-i})$ is not a transversal
and that $(b^{-1}a_{-i})\not\in C$.
Therefore the nerve $N$ of $W$ is equal to $px(a_{-i}^{-1}b)$ and it has greater syllable length than the nerve $N'=px$ of $W'$.
Take the core elements of $W'$, add to them $v_{i_n}=a_{-i}^{-1}b$ and declare the result the collection of core elements for $W$.
Proposition A follows now from the inductive hypothesis.

\smallskip
{\bf Subcase 4.D.} Suppose that $\hat z$ and $\hat u$ are nonempty.

Then $s_w=yz$, $s_t=yu$ and so
$h=pxa_{-i}^{-1}s_w^{-1}s_tbs_t^{-1}=pxa_{-i}^{-1}z^{-1}y^{-1}yubs_t^{-1}=pxa_{-i}^{-1}z^{-1}ubs_t^{-1}$.
Suppose $\hat y$ ends in a vertex of type $A_k$.
Let $z=f_1z_1$ and $u=f_2u_1$ where $f_1\in T_k$ is the label of the first edge of $\hat z$ and $f_2\in T_k$ is the label of the first edge of $u$.
Then clearly $f_1C\ne f_2C$ and so $f_1^{-1}f_2\not\in C$.
Notice also that $b\not\in C$ since if $b\in C$ then $s_tbs_t^{-1}$ stabilizes the last edge of $\hat u$ which contradicts the fact that $p$ is the reduced form for $h$ with respect to presentation (3).
Thus $W=pxa_{-i}^{-1}z_1^{-1}(f_1^{-1}f_2)u_1bs_t^{-1}$ is the normal form for $h$ with respect to presentation (1).
Since $b\not\in C$, Lemma 1.9(i) implies that $s_tb^{-1}$ is not a transversal. That is why the nerve $N$ of $W$ is equal to $pxa_{-i}^{-1}z_1^{-1}(f_1^{-1}f_2)u_1b$ and it has greater syllable length than the nerve $N'=px$ of $W'$.
Take the core elements of $W'$, add to them $v_{i_n}=b$ and declare the result the collection of core elements of $W$.
Proposition A follows now from the inductive hypothesis.
\smallskip
This completes the proof of Proposition A.
\enddemo

\proclaim {Corollary 1.11 (c.f. Corollary 4 from the Introduction)} Suppose $G=A_1\ast_C A_{-1}$ where the groups $G$ and $C$ are finitely generated. Suppose $H$ is a finitely generated subgroup of $G$ such that for any $g\in G$ we have $g^{-1}Hg\cap A_1= g^{-1}Hg\cap A_{-1}=\{1\}$. Then the subgroup $H$ is quasiisometrically embedded in $G$ (in particular, if $G$ is word hyperbolic then $H$ is quasiconvex in $G$).
\endproclaim
\demo {Proof}
Since $G$ and $C$ are finitely generated, the groups $A_1$ and $A_{-1}$ are also finitely generated. Fix a finite generating set ${\Cal C}$ of $C$ and a finite generating set $X_i$ containing ${\Cal C}$ of $A_i$ for $i=\pm 1$.
Put ${\Cal G}=X_1\cup X_{-1}$ to be the finite generating set of $G$.

Let $T, Y, Y_1$ and ${\Bbb B}$ be as in Proposition A.
Then $H$ is a free group on ${\Cal H}=E^{+}(B-Y_1)$ since $A_v=\{1\}$ for each $v\in VY_1$.
Suppose $h\in H$ and $U=U_1\dots U_2$ is a freely reduced word over ${\Cal H}=E^{+}(B-Y_1)$, $U_i\in {\Cal H}^{\pm 1}$. By Proposition A there is a reduced form
$W=v_1\dots v_m$ of $h$ with respect to the presentation $G=A_1\ast_C A_{-1}$ such that
$n\le m$.
On the other hand $m$ is the syllable length of $h$ with respect to the presentation $G=A_1\ast_C A_{-1}$. Therefore $l_{\Cal G}(h)\ge m$.

Thus $l_{\Cal H}(h)=n\le m\le l_{\Cal G}(h)$ and so $H$ is quasiisometrically embedded in $G$. 
\enddemo

\head 2. Word metric on fundamental groups of graphs of groups \endhead

\subhead Some auxiliarily facts \endsubhead

\proclaim {Lemma 2.1} Let $G$ be a word hyperbolic group generated by a finite set ${\Cal G}$. Let $w=w_1\dots w_t$ be a $K$-quasigeodesic word over ${\Cal G}$ where all $w_i$ are nonempty. Suppose for each $i=1,\dots ,t$ the word $u_i$ represents $\overline {w_i}$ and is $\lambda$-quasigeodesic.
Then for some constant $K'>0$ depending only on $K, \lambda$ the word
$w'=u_1\dots u_t$ is $K'$-quasigeodesic.
\endproclaim

\demo{Proof}
The statement of Lemma 2.1 is rather transparent and its proof is a standard exercise on quasiconvexity. Nevertheless the fact is of importance here and we will give a detailed argument.

Let $K_1=max (K,\lambda)$ and suppose any two $K_1$-quasigeodesics with the same endpoints in the Cayley graph of $G$ are $\epsilon$-Hausdorff-close.
Let $u$ be a subword of $w'$.
There are two possibilities.

\proclaim {Case 1} There is $u_i$ such that $u$ is a subword of $u_i$.
\endproclaim 
In this case, obviously, 
$$l(u)\le K_1\cdot l_{\Cal G}(\overline u) +K_1$$

\proclaim {Case 2} The word $u$ has the form $u=u_i'u_{i+1}\dots u_{j-1}u_j'$ where $i<j$, $u_i'$ is a terminal segment (perhaps empty) of $u_i$ and $u_j''$ is an initial segment (perhaps empty) of $u_j$.
\endproclaim

We want to show that for some constant $K'$

$$l(u_i'u_{i+1}\dots u_{j-1}u_j'')\le K' l_{\Cal G}(\overline {u_i'u_{i+1}\dots u_{j-1}u_j''})+K'\eqno (\dag)$$

There is a terminal segment $w_i'$ of $w_i$ and an initial segment $w_j''$ of $w_j$ such that $l_{\Cal G}(\overline {w_i'{u_i'}^{-1}})\le \epsilon$ and $l_{\Cal G}(\overline {{w_j''}^{-1}u_j''})\le \epsilon$.

Therefore $l_{\Cal G}(\overline {w_i'w_{i+1}\dots w_{j-1}w_j''})\le l_{\Cal G}(\overline {u_i'u_{i+1}\dots u_{j-1}u_j''})+2\epsilon$.  
We have 
$$l(u_i')\le K_1l_{\Cal G}(\overline {u_i'})+K_1\le K_1 (l_{\Cal G}(\overline {w_i'})+\epsilon)+K_1\le K_1(l({w_i'})+\epsilon)+K_1,$$

$$l(u_j'')\le K_1l_{\Cal G}(\overline {u_j''})+K_1\le K_1 (l_{\Cal G}(\overline {w_j''})+\epsilon)+K_1\le K_1(l({w_j''})+\epsilon)+K_1,$$

$$l(u_k)\le K_1l_{\Cal G}(\overline {u_k})+K_1\le K_1l({w_k})+K_1$$
and 
$$j-i\le l(w_i'w_{i+1}\dots w_{j-1}w_j'')+2.$$
Therefore
$$\eqalign {&l(u_i'u_{i+1}\dots u_{j-1}u_j'')\le K_1l(w_i'w_{i+1}\dots w_{j-1}w_j'')+K_1(i-j)+2(K_1+\epsilon)\le \cr 
&(K_1+1)l(w_i'w_{i+1}\dots w_{j-1}w_j'')+2(2K_1+\epsilon)}.$$
Put $K_2=2(2K_1+\epsilon)+1$. Then 
$$\eqalign{&l(u_i'u_{i+1}\dots u_{j-1}u_j'')\le K_2l(w_i'w_{i+1}\dots w_{j-1}w_j'')+K_2\le \cr &K_2K_1l_{\Cal G}(\overline {w_i'w_{i+1}\dots w_{j-1}w_j''})+K_2K_1+K_2\le \cr &K_2K_1(l_{\Cal G}(\overline {u_i'u_{i+1}\dots u_{j-1}u_j''})+2\epsilon)+K_2K_1+K_2\le \cr  &K'l_{\Cal G}(\overline {u_i'u_{i+1}\dots u_{j-1}u_j''})+K'}$$ where $K'=2K_2K_1+2K_2K_1\epsilon+K_2K_1+K_2$. Thus $(\dag)$ is established and Lemma 2.1 is proved.
\enddemo

\proclaim {Lemma 2.2} Suppose $G$ is a word hyperbolic group generated by a finite set ${\Cal G}$. Suppose $C_1,C_2\le G$ are virtually cyclic subgroups of $G$ such that is $C_1\cap C_2$ is finite. Let ${\Cal C}_i$ be a finite generating set of $C_i$. Assume that ${\Cal C}_i\subseteq {\Cal G}$.

Then 
\roster
\item there is a constant $\lambda >0$ such that whenever $U$ is a $d_{\Cal G}$-geodesic word such that $\overline U$ is shortest in the coset class $\overline U\cdot C_1$ and $V$ is a $d_{{\Cal C}_1}$-geodesic word, then the word $UV$ is $\lambda$-quasigeodesic with respect to $d_{\Cal G}$;
\item there is a constant $K>0$ such that whenever $u\in G$ is shortest in the double coset class $C_1uC_1$ and $c_1\in C$, the element $c_1u$ is at most $K$-away from any shortest element in the coset class $c_1uC_1$;
\item there is a constant $\lambda_1>0$ such that for any $d_{{\Cal C}_1}$-geodesic words $V, V'$ and any $d_{\Cal G}$-geodesic word such that $\overline U$ is shortest in the double coset class $C_1\overline U\cdot C_1$, then the word $VUV'$ is $\lambda_1$-quasigeodesic with respect to $d_{\Cal G}$;

\item there is a constant $K_1>0$ such that whenever $u\in G$ is shortest in the double coset class $C_1uC_1$ and $c_1\in C_1$, the word $c_1u$ is at most $K_1$-away from any shortest element in the coset class $c_1uC_1$;

\item there is a constant $K_2>0$ such that for any  $c_1\in C_1$ there is an element $c_2\in C_2$ with $l_{\Cal G}(c_2)\le K_2$ such that the element $u=c_1c_2$ is shortest in the coset class $c_1C_2$;
\item there is a constant $K_3>0$ such that for any $c_1\in C_1$  we have
$l_{\Cal G}(u)\ge l_{\Cal G}(c_1)/K_3-K_3$ where $u$ is shortest in the coset class $c_1C_2$;
\item there is a constant $\lambda_2>0$ such that for any $c_1\in C_1$ and $c_2\in C_2$ the word $V_1V_2$ is $\lambda_2$-quasigeodesic in the Cayley graph $\Gamma(G,{\Cal G})$ of $G$ where $V_i$ is a $d_{\Cal G}$-geodesic representative of $c_i$, $i=1,2$.

\endroster
\endproclaim
\demo {Proof}

(1), (2), (3) and (4) follow from the proof of Theorem C in [3]

(5) Let $y\in C_1$ be such that the cyclic group $<y>$ has finite index in $C_1$. Similarly, let $x\in C_2$ be such that the cyclic group $<x>$ has finite index in $C_2$. Fix a finite subset $T_1\subseteq C_1$ such that $C_1=T_1<y>$ and a finite subset $T_2\subseteq C_2$ such that $C_2=T_2<x>$.
 Observe that the statement (5) of Lemma 2.2 is obvious when at least one of the groups $C_1$, $C_2$ is finite.
From now on assume that they are both infinite. Thus $x,y$ are of infinite order and no nonzero power of $x$ is equal to a nonzero power of $y$ since $C_1\cap C_2$ is finite.
Let $c_1\in C_1$ be an arbitrary element. Then $c_1=t_1y^n$ for some $t_1\in T_1$.

Let $t_1y^n=uc_2$  where $c_2\in C_2$ and $u$ is shortest in the coset class $c_1C_2=t_1y^nC_2$. Thus $t_1y^n=ut_2x^k$ for some $t_2\in T_2$. Since $<x>$ is infinite and quasiconvex in $G$, there is a constant $K_2>0$ independent of $k$ and $n$ such that $u_1=ut_2$ is $K_2$-close to a shortest element in $ut_2<x>$.
Indeed, assume $ut_2=u_1=u'x^p$ where $u'$ is shortest in $ut_2<x>=u'<x>$.
It follows from the proof of Theorem C in [3] that there is a constant $N>0$ independent of $n,k,p$ such that $l_{Cal G}(u')+l_{\Cal G}(x^p)\le l_{\Cal G}(u_1)+N$. Suppose $p$ is such that $l_{\Cal G}(x^p)> l_{\Cal G}(t_2)+N$.
Then $l_{\Cal G}(u)\ge  l_{\Cal G}(u_1)-l_{\Cal G}(t_2)\le l_{Cal G}(u')+l_{\Cal G}(x^p)-N-l_{\Cal G}(t_2)>l_{Cal G}(u')$. Notice also that $uC_2=u'C_2$ which contradicts our choice of $u$. Thus $l_{\Cal G}(x^p)\le K_2=l_{\Cal G}(t_2)+N$.
 
We want to show that $|k|$ is small.
Let $Q_1,Y, U_1$ and $X$ be ${\Cal G}$-words representing $t_1$,$y$, $u_1$ and $x$.
It follows from (1) that there is $\lambda >0$ independent of $n$ and $k$ such that $UX^k$ and $Q_1Y^n$ are $\lambda$-quasigeodesics in the Cayley graph $\Gamma(G,{\Cal G})$ of $G$. Thus the paths $Q_1Y^n$ and $UX^k$ are $\epsilon$-hausdorff-close for some constant $\epsilon>0$.
If $|k|$ is greater than the number of elements in $G$ of length at most $\epsilon +2$ then there are numbers $n_1,n_2,k_1,k_2\ne 0$ such that
$zy^{n_2}z^{-1}=x^{k_2}$ where $z=x^{-k_1}y^{n_1}$.
Therefore 
$x^{-k_1}y^{n_1}y^{n_2}y^{-n_1}x^{k_1}=x^{k_2}$ and $y^{n_2}=x^{k_2}$.
This contradicts the fact that $<x>\cap <y>=\{1\}$.
Thus $|k|$ is bounded by a constant independent of $n$ which implies
statement (5) of Lemma 2.2.

(6) follows from (1).
\enddemo

\subhead Word metric on fundamental groups of graphs of groups \endsubhead

Suppose a word hyperbolic group $G$ is the fundamental group of a finite graph of groups ${\Bbb A}$ with respect to a maximum subtree $T$.

$$G=\pi_1({\Bbb A}, T) \eqno (3).$$

Assume that all edge groups $A_e$ are virtually cyclic.

Then $G$ has a presentation
$$G=\underset{v\in VA}\to {(\ast A_v)}\ast F(E^{+}A)/\{e=1, e\in ET; \alpha_e(a)e=e\omega_e(a), e\in E^{+}A, a\in A_e\}.$$ 

For each $e\in EA$ we fix $x_e\in A_e$ such that $<x_e>$ has finite index in $A_e$.
Denote $x_{e,\alpha}=\alpha_e(x_e)\in A_{\partial_0(e)}$ and
$x_{e,\omega}=\omega_e(x_e)\in A_{\partial_1(e)}$.

For each vertex $v\in VA$ we fix a finite generating set $Z_v$ closed under taking inverses. We may assume that for each edge $e$ of ${\Bbb A}$ originating from $v$ the set  $Z_v$ includes the generator $c_{e,\alpha}$ of the subgroup of finite index in $\alpha_e(A_e)$.
Put 
$$Z=\underset{v\in VA}\to\cup Z_v\bigcup \{e| e\in E(A-T)\}$$ 
Then $Z$ is a finite generating set for $A$.
Put $Z'=Z\cup ET$. It is another finite generating set for $G$ (every $e\in ET$ represents the trivial element of $G$).
Any $d_{Z'}$-geodesic word $w$ contains no letters $e\in ET$ and so it is a word over $Z$. Clearly it is $d_Z$-geodesic.
Thus $d_Z$ and $d_{Z'}$ coincide on $G$.

\proclaim {Lemma 2.3} Suppose $W=UeV$ where $e\in EA$, $v_0=\partial_0(e)$, $v_1=\partial_1(e)$, $U$ is a $d_{Z_{v_0}}$-geodesic word, $v_1=\partial_1(e)$, $V$ is a $d_{Z_{v_0}}$-geodesic word.
Suppose $\overline V=\in \omega_{e}(A_{e})$.

Assume that $W$ is a $K$-quasigeodesic in the $d_{Z'}$-metric for some $K>0$.
Then there is $K'>0$ independent of $\overline U$, $\overline V$, such that $W'=U_1e$ is $K'$-quasigeodesic where $U_1$ is a $d_{Z_{v_0}}$-geodesic word representing $\overline {U}\alpha_e(\omega_e)^{-1}(\overline V)$.
 
\endproclaim 

\demo {Proof}
The subgroup $A_{\partial_0(e)}=A_{v_0}$ is e quasiconvex in $G$. Thus there is a constant $K_1$ such that $ye$ is $K_1$-quasigeodesic in the $d_{Z'}$-metric for any $d_{Z_{v_0}}$-geodesic word $y$.
Let $V_1$ be a $d_{Z_{v_0}}$-geodesic word representing $\alpha_e(\omega_e)^{-1}(\overline V)$. Then $\overline {V_1e}=\overline {eV}$. 
Therefore by Lemma 2.1 the word $W_2=UV_1e$ is $K_2$-quasigeodesic in the $d_{Z'}$-metric for some constant $K_2$.

Recall that $U_1$ is a $d_{G_{v_0}}$-geodesic representative of $\overline {UV_1}$. Since $A_{G_{v_0}}$ is quasiconvex in $G$, we know that $U_1$ is a $K_3$-quasigeodesic in the Cayley graph of $G$. Thus by Lemma 2.1 $W'=U_1e$ is a $K'$-quasigeodesic in the $d_{Z'}$-metric for some constant $K'$.
\enddemo

\proclaim {Lemma 2.4} Let $F$ be the subgroup of $G$ generated by $EA$ that is $F$ is a free group on $E^{+}(A-T)$. Then
\roster
\item "(a)" The subgroup $F$ is isometrically embedded in $G$ that is any freely reduced word in $E^{+}(A-T)$ is $d_Z$-geodesic.
\item "(b)" For each $K>0$ and $M>0$ there is $K_1>0$ such that whenever $W=W_1\dots W_s$ is a $K$-quasigeodesic in the $d_{Z'}$-metric and $U_0$,\dots ,$U_{s}$ are words in $E^{+}T$ of length at most $M$, the word 
$$W'=U_0W_1U_1\dots U_{s-1}W_sU_s$$ is $K_1$-quasigeodesic in the $d_{Z'}$-metric.

\endroster
\endproclaim

\demo {Proof}
Statement (a) follows obviously from the properties of HNN-extensions.
\smallskip
Statement (b) is obvious. 
\enddemo

\proclaim {Proposition B} There is a constant $K>0$ such that for any $g\in G$ there is a $K$-quasigeodesic with respect to $d_{Z'}$ word $W$ representing $g$ of the form
$$W=W_1\dots W_n$$ where each $W_k$ is either $e^{\pm 1}$ for some $e\in EA$ or $W_k$ is a $d_{Z_v}$-geodesic word for some $v\in VA$ and 
$$\overline {W_1},\dots ,\overline {W_n}$$ is a reduced form for $h$ with respect to presentation (3).
\endproclaim

\demo {Proof}
Let $W$ be a $Z$-geodesic word representing $g$.
We will transform $W$ to the required form in several steps.
\smallskip
{\bf Step 1} We can write $W$ as $W=Q_1\dots Q_m$ where for $k=1,\dots ,m$ each $Q_k$ is either an edge of $A-T$ or it represents an element of a vertex group of $A$ and whenever $1\le i<j\le m$, the word $Q_i\dots Q_j$ does not represent an element of a vertex group of $A$.
Notice that we do not claim that each $Q_i$ is a word in generators of some vertex group.
Now each $A_v$ is quasiconvex in $G$ since all the edge groups are virtually cyclic [11].
Let $K_1>0$ be such that for any $v\in VA$ any $d_{Z_v}$ geodesic word is $K_1$-quasigeodesic in the $d_Z$-metric.
For each $k=1,\dots ,m$ such that $\overline {Q_k}\in A_v$ we find a $d_{Z_v}$-geodesic representative $U_k$ of $\overline {Q_k}$.
For other $Q_k$ we put $U_k=Q_k$.
Let $W_1=U_1\dots U_m$. Then by Lemma 2.1 the word $W_1$ is $K_2$-quasigeodesic in the $d_Z$-metric for some constant $K_2$ independent of $g$.

\smallskip
{\bf Step 2} Now between every $U_k, U_{k+1}$ representing elements of vertex groups $v_{k}$ and $v_{k+1}$ of $A$ we insert the reduced edge-path $r_k$ in $T$ from $v_k$ to $v_{k+1}$.
Between every $U_k, U_{k+1}$ such that $\overline {U_k}\in A_{v_k}$ and $U_{k+1}=e\in E(A-T)$ we insert the reduced edge-path $r_k$ in $T$ from $v_k$ to the initial vertex of $e$.
Between every $U_k, U_{k+1}$ such that $U_{k}=e\in E(A-T)$ and $\overline {U_{k+1}}\in A_{v_{k+1}}$ we insert the reduced edge-path $r_k$ in $T$ from the terminal vertex of $e$ to $v_k$.
We put $r_0$ to be the reduced edge-path from $d_1$ to the initial vertex of $U_1$ when $U_1$ is an edge of $A-T$ and we put $r_0$ to be the reduced edge-path from $d_1$ to the vertex $v_1$ when $\overline {U_k}\in A_{v_1}$.
Analogously, we put $r_m$ to be the reduced edge-path from the terminal vertex of $U_m$ to $d_1$ when $U_m$ is an edge of $A-T$ and we put $r_m$ to be the reduced edge-path from $v_m$ to $d_1$ when $\overline {U_m}\in A_{v_m}$
Then $$r=r_0,\overline {U_1},r_1,\dots r_{m-1},\overline {U_m},r_m$$ is a loop at $d_1$ in the  graph of groups ${\Cal A}$ which represents $g$.
Observe that each edge-path $r_i$ has length at most $N_0$ where $N_0$ is the number of oriented edges of $T$ since $T$ is a tree and $r_i$ has no backtrackings. Therefore  by Lemma 2.4 the word 
$$W_2=r_0U_1r_2\dots  r_{m-1}U_mr_m$$ is a $K_3$-quasigeodesic with respect to $d_{Z'}$-metric where $K_3$ is some constant independent of $g$.

\smallskip
{\bf Step 3}
For each $k=1,\dots ,m$ we find the maximal initial segment $r_k'$ of $r_k$ such that ${r_k'}^{-1}=r_{k-1}''$ is a terminal segment of $r_{k-1}$ and ${r_k'}^{-1}U_kr_k$ represents an element $u_k'$ of a vertex group of $A$. Notice that $r_k'$ and $r_k''$ are disjoint subwords of $r_k$ since otherwise there is a subword $U_i\dots U_j$ of $W_1$, $i<j$, representing an element of a vertex group.
 
 Replace in $r$ each loop ${r_k'}^{-1},\overline {U_k},r_k'$  by $u_k'$.

This gives us the path
$$r'=r_0'u_1'r_2'\dots r_{m-1}'u_m'r_0'$$ such that 
\roster
\item each $u_k'$ is either an edge of $A-T$ or a nontrivial element of a vertex group of ${\Cal A}$,
\item whenever $i<j$ the element $u_i'\dots u_j'$ is not in a vertex group of ${\Cal A}$
\item whenever the last edge $e$ of $r_i'$ is inverse to the first edge of $r_{i+1}'$, then $\overline {u_i'}\not\in \omega_e(A_e)$.
\endroster

Let $U_i'=e$ whenever $u_i'=e\in E(A-T)$ and let $U_i'$ be a $d_{Z_v}$-geodesic representative of $u_i'$ whenever $u_i'\in A_v$. For each $i=1,\dots ,m$ replace the subword $r_{k-1}''U_kr_k'$ of $W_2$ by the word $U_k'$. This produces a new word
$$W_3=r_0'U_1'r_2'\dots r_{m-1}'U_m'r_0'.$$
Since the vertex groups are quasiconvex in $G$, Lemma 2.1 implies that
$W_3$ is a $K_3$-quasigeodesic in the $d_{Z'}$-metric where $K_3$ is a constant independent of $g$.
Notice that $r'$ is close to being a {\it normal form} of $g$ with respect to presentation (3). The only conditions of Definition 1.1 which are possibly not satisfied are conditions 5) and 6). 
\smallskip
{\bf Step 4} 
The sequence $r'$ can be broken into maximal pieces
$r'=P_1\dots P_{s'}$ where each $P_k$ has the form

$$\eqalign {&g_{i_k},e_{i_k},\omega_{e_{i_k}}(c_{i_k})\alpha_{e_{i_k+1}}(c_{i_k+1}^{-1}),e_{{i_k}+1},\omega_{e_{i_k+1}}(c_{i_k+1})\alpha_{e_{{i_k}+2}}(c_{i_k+2}^{-1}),e_{{i_k}+2},\dots \cr
&\dots ,e_{{j_k}-1},\omega_{e_{j_k-1}}(c_{j_k-1})\alpha_{e_{j_k}}(c_{j_k}^{-1}), e_{j_k}, \omega_{e_{j_k}}(c_{j_k}) } \eqno (5)$$

where $c_s\in A_{e_s}$ for $s=i_k,\dots,j_k$, $g_{i_k}\in A_{\partial_0(e_{i_k})}$ and $g_{i_k}\alpha_{e_{i_k}}(c_{i_k})$ cannot be "pulled to the left" that is either $i_k=1$ and $e_{i_k}=e_1$ is the first edge of $r'$ or $i_k>1$ and $g_{i_k}\alpha_{e_{i_k}}(c_{i_k})\not\in \omega_{e_{i_k-1}}(A_{e_{i_k-1}})$.
To each $P_k$ there is a corresponding subword $\hat P_k$ of $W_2$ of the form
$$\hat P_k=V_{i_k}e_{i_k}V_{i_k+1}e_{i_k+1}\dots e_{j_k}V_{j
_k+1} \eqno (6)$$

Observe that for each $k$ the path $e_{i_k},\dots ,e_{j_k}$ has no backtrackings.
Indeed, suppose $e_{s+1}=e_s^{-1}$. Then by construction of $r'$ we have  $1\ne\omega_{e_{s}}(c_{s})\alpha_{e_{s+1}}(c_{s+1}^{-1})=\overline {U_t}$ for some $t$.  
On the other hand, equation (5) implies that 
$\omega_{e_{s}}(c_{s})\alpha_{e_{s+1}}(c_{s+1}^{-1})=\omega_{e_{s}}(c_{s})\omega_{e_{s}}(c_{s+1}^{-1})=\omega_{e_{s}}(c_{s}c_{s+1}^{-1})$ and therefore $e_s\omega_{e_{s}}(c_{s})\alpha_{e_{s+1}}(c_{s+1}^{-1})e_{s+1}=\alpha_{e_{s}}(c_{s}c_{s+1}^{-1})$. This contradicts property 3) of $r'$.

For each $k$ let $l_k$ be the maximal among $\{i_k,\dots ,j_k-1\}$ such that for $s=i_k,\dots ,l_k-1$ the groups $\omega_{e_{s}}(A_{e_{s}})$ and $\alpha_{e_{s+1}}(A_{e_{s+1}})$ have infinite intersection. If there are no such indices among $\{i_k,\dots ,j_k-1\}$, we put $l_k=i_k$.

\enddemo
\proclaim {Claim 1}
We have $l_k-i_k< N+1$ where $N$ is the number of oriented edges in the graph $A$. 
\endproclaim

Indeed, suppose not and $l_k-i_k\ge N+1$
Then there is a subpath $e_{s_1},\dots ,e_{t_1}$ of $e_{i_k},\dots ,e_{l_k}$, $s_1\ne t_1$, such that $e_{s_1}=e_{t_1}$. We know that for each 
$s=i_k,\dots ,l_k$ the groups $<x_{s,\omega}>$ and $<x_{s+1,\alpha}>$ are infinite and commensurable (i.e. they are finite extension of a common infinite cyclic subgroup). Thus we can find $M\ne 0$ and $M_1\ne 0$ such that
$x_{s_1,\alpha}^{M}e_{s_1}\dots e_{t_1-1}=e_{s_1}\dots e_{t_1-1}x_{t_1,\alpha}^{M_1}=e_{s_1}\dots e_{t_1-1}x_{s_1,\alpha}^{M_1}$ in $G$.
Notice that the edge-path $e_{s_1},\dots ,e_{t_1}$ contains an edge which is not equal to the trivial element in $G$.
Indeed, if they are all trivial then, since $e_{s_1}=e_{t_1}$, the path  $e_{s_1},\dots ,e_{t_1}$ contains a backtracking $e,e^{-1}$ which is impossible.
Thus  $e_{s_1}\dots e_{t_1-1}$ represents a non-trivial element $f\in F$.
Therefore $x^Mf=fx^{M_1}$ where $x=x_{s_1,\alpha}=x_{t_1,\alpha}$.
Since $G$ is word hyperbolic, this implies $<f>\cap <x>\ne \{1\}$ which is impossible by standard properties of graphs of groups. 
Thus we have established that $l_k-i_k< N+1$ and Claim 1 is proved.
\smallskip
Observe that $P_k$
represents the same element of $G$ as
$$\eqalign {P_k'&=g_{i_k},e_{i_k},\omega_{e_{i_k}}(c_{i_k})\alpha_{e_{i_k+1}}(c_{i_k+1}^{-1}),e_{{i_k}+1},\omega_{e_{i_k+1}}(c_{i_k+1})\alpha_{e_{{i_k}+2}}(c_{i_k+2}^{-1}),e_{{i_k}+2},\dots \cr
&\dots ,e_{{l_k}-1},\omega_{e_{l_k-1}}(c_{l_k-1})\alpha_{e_{l_k}}(c_{l_k}^{-1}), e_{l_k}, \omega_{e_{l_k}}(c_{l_k}),e_{l_k+1},e_{l_k+2},\dots ,e_{j_k}}$$
By definition of $l_k$ the groups $\omega_{e_{l_k}}(A_{e_{l_k}})$ and $\alpha_{e_{l_k+1}}(A_{e_{l_k+1}})$ have finite intersection. 

Thus by Lemma 2.2 the word $V_{l_k}'V_{l_k}''$ is $\lambda$-quasigeodesic in the Cayley graph of $G$ where $V_{l_k}'$ is a $d_{Z_v}$-geodesic representative of $\omega_{e_{l_k}}(c_{l_k})$ and $V_{l_k}''$ is a $d_{Z_v}$-geodesic representative of $\alpha_{e_{l_k+1}}(c_{l_k}^{-1})$, $v=\partial_1(e_{l_k})=\partial_0(e_{l_k+1})$.

In the word $W_3=r_0'U_1'r_1'\dots r_{m-1}'U_m'r_m'$ for each $k=1,\dots ,s'$ we substitute the word $V_{l_k+1}$ representing the element $\omega_{e_{l_k}}(c_{l_k})\cdot \alpha_{e_{l_k+1}}(c_{l_k+1})$ by the word $V_{l_k}'V_{l_k}''$.
This produces a new word
$W_4=r_0''U_1''r_1''\dots r_{m-1}''U_m''r_m''$ which is $\lambda_1$-quasigeodesic by Lemma 2.1.

Besides, we know that the edge-path $e_{i_k},\dots ,e_{j_k}$ has no backtrackings and therefore the lengths of the segments, into which elements of $E(A-T)$ divide it, do not exceed $N$ where $N$ is the number of edges in $A$. Thus by Lemma 2.4 the word $e_{k_0+1},e_{k_0+2},\dots ,e_j$ is $\lambda_4$-quasigeodesic in $d_{Z'}$-metric where $\lambda_4$ is some constant independent of $g$. 

For each $k=1,\dots ,s'$ let $V_{i_k}'$ be a $d_{Z_{\partial_0(e_i)}}$-geodesic representative of $g_{i_k}\alpha_{e_{i_k}}(c_{i_k})$.
Now for each $k=1,\dots ,s'$ we replace the subword $V_{l_k}''e_{l_k+1}V_{l_k+1}e_{l_k+2}V_{l_k+2}\dots V_{j_k}e_{j_k}V_{j_k+1}$ of $W_4$
by $e_{l_k+1}e_{l_k+2}\dots e_{j_k}$ and the subword $V_{i_k}e_{i_k}V_{i_k+1}e_{{i_k}+1}V_{i_k+2}\dots V_{l_k}e_{l_k}V_{l_k}'$ of $W_4$ by
$V_{i_k}'e_{i_k}e_{i_k+1}\dots e_{l_k}$ to get a word $W_4$. Since $l_k-i_k\le N+1$, Lemma 2.1 implies that the new word $W_4$ is $\lambda_5$-quasigeodesic in the $d_{Z'}$-metric where $\lambda_5$ is some constant independent of $g$.
It also follows from the construction that $W_4$ satisfies all requirements of Definition 1.1 except, possibly, condition 6).
\smallskip
{\bf Step 5}
Let $r^{(3)}$ be the path in ${\Bbb A}$ corresponding to the word $W_4$, that is $r^{(3)}$ is obtained from $W_4$ by "barring" every letter. 

Then $r^{(3)}$ can be broken into maximal pieces 

$r^{(3)}=R_1\dots R_{s''}$ where each piece $R_k$ has the form

$$g_{l_0(k)}, r_{1,k}, g_{l_1(k)},r_{2,k}\dots r_{j,k} g_{l_j(k)}$$
where each $g_{l_i(k)}\in A_{v_{i(k)}}-\{1\}$ is a nontrivial vertex group element, each $r_{i,k}$ is an edge-path with trivial vertex group elements inserted between the consecutive edges and for each $i=0,\dots ,j-1$ the element $g_{l_i(k)}$ can be "pulled through" the edge-path $r_{i+1,k}r_{i+2,k}\dots r_{j,k}$ to the element of $A_{v_{j(k)}}$.

Recall that each edge-path $r_{i,k}$ is without backtracks by construction of $W_4$. Moreover $r_{1,k}r_{2,k}\dots r_{j,k}$ also does not have backtracks. Indeed, if the last edge $e$ of $r_{i,k}$ is inverse to the first edge of $r_{i+1,k}$ then the element $g_{l_i(k)}$ can be "pulled to the left" trough the edge $e$ which contradicts the properties of $W_4$.

Find minimal $i$ (if any) such that $g_{l_i(k)}$ is of infinite order and denote it $i_0$. If there are no such $i$, put $i_0=l_j(k)$.

{\bf Claim 2} The length of $r_{i_0+1,k}r_{i_0+2,k}\dots r_{j,k}$ is at most $N$ where $N$ is the number of oriented edges in the graph $A$.

The proof is exactly the same as that of Claim 1.

Notice that there is a uniform bound on the lengths of elements of finite order in vertex groups which come are images of elements of finite order in edge groups.
For each $R_k$ there is a corresponding subword $\hat R_k=V_{l_0(k)}r_{1,k} V_{l_1(k)}\dots V_{l_j(k)}$ of $W_4$.
Let $V_{l_j(k)}'$ be a $Z_{v_{j(k)}}$-geodesic word such that
$r_{1,k}r_{2,k}\dots r_{j,k}V_{l_j(k)}'$ represents the same element of $G$ as $R_k$. Notice that $\overline {V_{l_j(k)}'}=c\overline {V_{l_j(k)}}$ for some $c\in \omega_e(A_e)$ where $e$ is the last edge of $r_{j,k}$. Thus  $\overline {V_{l_j(k)}'}\not\in \omega_e(A_e)$ since $\overline {V_{l_j(k)}}\not\in \omega_e(A_e)$.

Substitute every $\hat R_k$ in $W_4$ by $r_{1,k}r_{2,k}\dots r_{j,k}V_{l_j(k)}'$ to get a new word $W_5$. Applying Lemma 2.3 at most $N$ times we conclude that $W_5$ is $\lambda_6$-quasigeodesic with respect to $d_{Z'}$ where $\lambda_6$ does not depend on $g$. It follows from the construction that $W_5$ corresponds to a reduced form of $g$ with respect to presentation (3). 
This completes the proof of Proposition B.

\head 3. Proofs of Theorem A and Theorem B\endhead

\proclaim {Theorem B} Let $$G=A_1\ast_C A_{-1} \eqno(7)$$ be a word hyperbolic group where the group $C$ is virtually cyclic (and therefore $A_1$ and $A_{-1}$ are quasiconvex in $G$). 
Let $H$ be a finitely generated subgroup of $G$.

Then $H$ is quasiconvex in $G$ if and only if for each $g\in G$ and $i=\pm 1$ the subgroup $g^{-1}Hg\cap A_i$ is quasiconvex in $A_i$.
\endproclaim

Before proceeding with the proof of Theorem B we need the following

\proclaim {Lemma 3.1} Let $G$, $A_1$, $A_{-1}$ and $C$ be as above.
Let $X_i$ be a finite generating set of $A_i$ containing a finite generating set ${\Cal C}$ of $C$. Put ${\Cal G}=X_1\cup X_{-1}$ to be the finite generating set for $G$. Then the following holds.
\roster
\item There exists $K>0$ such that each $g\in G$ has a $K$-quasigeodesic (with respect to $d_{\Cal G}$) representative of the form
$$u=u_0\dots u_s$$
where 
{\roster
\item "(i)" each $u_k$ is a $d_{X_j}$-geodesic word for some $j\in \{\pm 1\}$ such that
when $k>0$ the element $\overline {u_k}$ does not belong to $C$ and it is shortest (with respect to $d_{X_j}$) in the coset class $C\overline {u_k}$;
\item "(ii)" If $\overline {u_k}\in A_j-C$ then $\overline {u_{k-1}}\in A_{-j}-C$, $k=1,\dots ,s$.  
\item "(iii)" If $s=0$ and $g=c\in C$ then $w_0$ is a $d_{X_1}$-geodesic representative of $c$.
\endroster}

\item  Suppose $N>0$ is a fixed number. Then there is a constant $K_1>0$ such that the following holds.
Suppose $g\in G$ and $w=w_0\dots w_s$ is a word such that
{\roster
\item "(i)" $\overline w =g$;
\item "(ii)" each $w_k$ is a a $d_{X_j}$-geodesic word for some $j\in \{\pm 1\}$ such that for $k>0$ the element $\overline {u_k}$ does not belong to $C$ ;
\item "(iii)" If $\overline {w_k}\in A_j-C$ then $\overline {w_{k-1}}\in A_{-j}-C$, $k=1,\dots ,s$;
\item "(iv)" If $k>0$, $\overline {w_k}=c'v_k$, $\overline {w_{k-1}}=v_{k-1}''c''$ where $c',c''\in C$, $v_k'$ is shortest (with respect to $d_{X_j}$ in the coset class $C\overline {w_k}$ and $v_{k-1}''$ is shortest (with respect to $d_{X_{-j}}$) in the coset class $\overline {w_{k-1}}C$ then

$$ l_{\Cal C}(c''c')\ge l_{\Cal C}(c'')+l_{\Cal C}(c')-N$$
 \endroster}
\endroster

Then $l_{\Cal G}(g)\ge K_1 l(w)+K_1$.
\endproclaim 
\demo {Proof}

(1) This is a more or less immediate corollary of Proposition B applied to the group $G$.
If $g=c\in C$, that is $w$ is a $d_{\Cal C}$-geodesic representative of $c$, the statement of Lemma 3.1(1) is obvious.
Assume from now on that $\overline w\not\in C$.
By Proposition B there exists $K>0$ such that every element $g\in G$ has a $d_{\Cal G}$-quasigeodesic representative of the form

$w=w_0\dots w_s$ where

\roster
\item  each $w_k$ is a $d_{X_j}$-geodesic word for some $j\in \{\pm 1\}$ such that the element $\overline {w_k}$ does not belong to $C$;
\item  if $\overline {w_k}\in A_j-C$ then $\overline {w_{k-1}}\in A_{-j}-C$, $k=1,\dots ,s$. 
 
\endroster

We will transform $w$ to the required form in several steps.

{\bf Step 1.}

For each $k=0,\dots ,s$ express $\overline {w_k}\in A_j-C$ as $\overline {w_k}=\overline {x_kv_kz_k}$ where $x_k, z_k$ are $d_{\Cal C}$-geodesic words, $v_k$ is $d_{X_j}$-geodesic word and $\overline {v_k}$ is shortest in the double coset class $C\overline {w_k}C$.
It follows from Lemma 2.2(3) that there is a constant $K_1>0$ independent of $g$ such that each $x_kv_kz_k$ is a $K_1$-quasigeodesic in the $d_{X_j}$ metric. Since $A_1$ and $A_{-1}$ are quasiconvex in $G$, there is $K_2>0$ independent of $g$ such that $x_kv_kz_k$ is $K_2$-quasigeodesic in $d_{\Cal G}$-metric.
Replace every $w_k$ by $x_kv_kz_k$ to get a word

$$w'=x_0v_0z_0x_1v_1z_1x_2\dots x_{s-1}v_{s-1}z_{s-1}x_sv_sz_s.$$

By Lemma 2.1 there is $\lambda >0$ independent of $g$ such that $w'$ is $\lambda$-quasigeodesic in the $d_{\Cal G}$-metric.

{\bf Step 2.} For each $k=0,\dots ,s-1$ we find a $d_{\Cal C}$-geodesic word $y_k$ representing the element $\overline {z_kx_{k+1}}$. Since $C$ is quasiconvex in $G$, there is a constant $K_2>0$ such that each $y_k$ is $K_2$-quasigeodesic with respect to $d_{\Cal G}$.
Replace for each $k=0,\dots ,s-1$ the word $z_kx_{k+1}$ by $y_k$ to get a new word
$$w''=x_0v_0y_0v_1y_1v_2\dots y_{s-2}v_{s-1}y_{s-1}v_sz_s.$$
By Lemma 2.1 there is $\lambda_1 >0$ independent of $g$ such that $w''$ is $\lambda_1$-quasigeodesic in the $d_{\Cal G}$-metric.

{\bf Step 3.} Express $\overline {v_sz_s}\in A_j-C$ as $\overline {v_sz_s}=\overline {q_su_s}$ where $u_s$ is a $d_{X_j}$-geodesic word, $q_s$ is a $d_{\Cal C}$-geodesic word and $\overline {q_s}$ is shortest (with respect to $d_{X_j}$) in the coset class $C\overline {v_sz_s}$.
Then express $\overline {v_{s-1}y_{s-1}q_s}\in A_{-j}-C$ as $\overline {v_{s-1}y_{s-1}q_s}=\overline {q_{s-1}u_{s-1}}$ where $u_{s-1}$ is a $d_{X_{-j}}$-geodesic word, $q_{s-1}$ is a $d_{\Cal C}$-geodesic word and $\overline {q_{s-1}}$ is shortest (with respect to $d_{X_{-j}}$) in the coset class $C\overline {v_{s-1}y_{s-1}q_s}$.
And so on. Finally, express $\overline {x_0v_0y_0q_1}\in A_i-C$
as $\overline {x_0v_0y_0q_1}=\overline {u_0}$ where $u_0$ is a $d_{X_i}$-geodesic word.

By Lemma 2.2(2) there is a constant $K_3>0$ independent of $g$ such that $l_{\Cal C}(\overline {q_k})\le K_3$, $k=1,\dots ,s$. Between each $y_{k-1}$ and $v_k$ is $w''$ we insert a word $q_kq_k^{-1}$ to get a word

$$w'''=x_0v_0q_1q_1^{-1}v_1y_1q_2q_2^{-1}\dots v_{s-1}y_{s-1}q_sq_s^{-1}v_sz_s.$$

Since $l_{\Cal C}(\overline {q_k})\le K_3$, the word $w'''$ is $\lambda_2$-quasigeodesic with respect to $d_{\Cal G}$ where $\lambda_2$ is a constant independent of $g$.

{\bf Step 4.} Recall that $A_1$ and $A_{-1}$ are quasiconvex in $G$.
Therefore there is a constant $K_4>0$ such that any $d_{X_i}$-geodesic word is $K_4$-quasigeodesic with respect to $d_{\Cal G}$. So there is $K_5>0$ independent of $g$ such that the word and $q_{k}u_k$ is $K_5$-quasigeodesic with respect to $d_{\Cal G}$ for $k=1,\dots ,s$.

Replace in $w'''$ each $v_{k}y_{k}q_{k+1}$ by $q_{k}u_{k}$ for $k=1,\dots ,s-1$, replace $x_0v_0q_1$ by $u_0$ and replace $v_sz_s$ by $q_su_s$ to get the word
$$w^{(4)}= u_0q_1^{-1}q_1u_1q_2^{-1}q_2u_2\dots q_{s-1}u_{s-1}q_s^{-1}q_su_s.$$
By Lemma 2.1 $w^{(4)}$ is $\lambda_3$-quasigeodesic for some constant $\lambda_3>0$ independent of $g$.

{\bf Step 5.} Finally we replace each $q_k^{-1}q_k$ by the empty word to get the word

$$u=u_0u_1\dots u_s.$$ 
By Lemma 2.1 $w^{(5)}$ is $\lambda_4$-quasigeodesic for some constant $\lambda_4>0$ independent of $g$.
It follows from the construction that $\overline u=g$ and that $u$ satisfies all the requirements of Lemma 3.1 (1).
This completes the proof of part (1) of Lemma 3.1.

(2) Let $w=w_0\dots w_s$ be as in Lemma 3.1 (2).
We will show that it can be transformed to a quasigeodesic $u=u_0\dots u_s$ as in Lemma 3.1(1) without loosing too much length.
Assume $g=\overline w\not\in C$.

For each $k=0,\dots ,s$ express $\overline {w_k}\in A_j-C$ as $\overline {x_kv_kz_k}$ where $x_k,z_k$ are $d_{\Cal C}$-geodesic words, $v_k$ is a $d_{X_j}$-geodesic word such that $\overline {v_k}$ is shortest in the double coset class $C\overline {w_k}C$.
Lemma 2.2 implies that $x_kv_kz_k$ is $\lambda$-quasigeodesic with respect to $d_{X_j}$ where $\lambda>0$ does not depend on $g$.
Obviously, $l(w_k)\le l(x_kv_kz_k)$ since $w_k$ is $d_{X_j}$-geodesic.
Put 
$$w'=x_0v_0z_0x_1v_1z_1\dots x_sv_sz_s.$$
Then $l(w)\le l(w')$.

 Notice that by Lemma 2.2 $\overline {x_kv_k}$ is at most $K_1$-away from the shortest element in the coset class $\overline {w_k}C$ and ${\overline {v_kz_k}}^{-1}$ is at most $K_1$ away from the inverse of the shortest element in $C\overline {w_k}$ (here $K_1>0$ is a constant independent of $g$).
Therefore there is $K_2>0$, depending on $N$ but independent of $g$, such that

$$l_{\Cal C}(\overline {z_{k-1}x_k})\ge l_{\Cal C}(\overline{z_{k-1}})+l_{\Cal C}(\overline{x_k})-K_2.$$

Take $y_k$ to be a $d_{\Cal C}$-geodesic representative of $\overline {z_{k-1}x_k}$. Then $l(z_{k-1}x_k)\le l(y_k)+K_2$.
Replace each $z_{k-1}x_k$ by $y_k$ in $w'$ to get 
$$w''=x_0v_0y_1v_1y_2\dots y_sv_sz_s.$$

Then $l(w')\le l(w'')+sK_2\le l(w'')+l(w'')K_2=(K_2+1)l(w'')$ and therefore
$l(w)\le l(w')\le (K_2+1)l(w'')$.
Express $\overline {v_sz_s}\in A_j-C$ as $\overline {q_su_s}$ where $q_s$ is $d_{\Cal C}$-geodesic, $u_s$ is $X_j$-geodesic and $\overline {u_s}$ is shortest with respect to $d_{X_j}$ in the coset class $C\overline {v_sz_s}$.
Then express $\overline {v_{s-1}y_sq_s}$ as $\overline {q_{s-1}u_{s-1}}$ where $q_{s-1}$ is $d_{\Cal C}$-geodesic, $u_{s-1}$ is $X_j$-geodesic and $\overline {u_{s-1}}$ is shortest with respect to $d_{X_{-j}}$ in the coset class $C\overline {v_{s-1}y_sq_s}$. And so on.
Finally, we rewrite $\overline {x_0v_0y_1q_1}\in A_i-C$ as $\overline {u_0}$ where $u_0$ is $d_{X_i}$-geodesic.
Recall that $l_{\Cal C}(\overline {q_k})\le K_1$. 
Therefore there is $K_3>0$ independent of $g$ such that

$$l(v_{k}y_{k+1})\le K_3 l(u_k), \ \ k=1,\dots ,s-1$$
$$l(v_sz_s)\le K_3 l(u_s)$$
$$l(x_0v_0y_1) \le K_3 l(u_0).$$

Put $u=u_0\dots u_s$. Then $l(w'')\le K_3l(u)$ and therefore $l(w)\le K_3(K_2+1)l(u)$. It is clear from the construction that $\overline u=g$ and that $u$ satisfies the requirements of Lemma 3.1(1) and therefore it is $K$-quasigeodesic with respect to $d_{\Cal G}$.
Thus $l(w)\le K_3(K_2+1)l(u)\le KK_3(K_2+1)l_{\Cal G}(g)+ KK_3(K_2+1)$ 
This completes the proof of Lemma 3.1.
\enddemo
 
\demo {Proof of Theorem B}
\enddemo

Suppose $H$ is quasiconvex in $G$. Then for every $g\in G$ and $i=\pm 1$ the subgroups $g^{-1}Hg$ and $A_i$ are quasiconvex in $G$. Therefore their intersection $g^{-1}Hg\cap A_i$ is quasiconvex in $G$. Since $A_i$ is quasiconvex in $G$ and $g^{-1}Hg\cap A_i\le A_i$ this implies that $g^{-1}Hg\cap A_i$ is quasiconvex in $A_i$.

From now on we assume that $H$ is a finitely generated subgroup of $G$ such that for any $g\in G$ and $i=\pm 1$ the subgroup $g^{-1}Hg\cap A_i$ is quasiconvex in $A_i$.
If $H$ is conjugate in $G$ to a subgroup of $A_i$ then $H$ is quasiconvex in a conjugate of $A_i$ and so in $G$. Thus we may assume that $H$ is not conjugate to a subgroup of $A_i$.

We recall some notations from section 1 which will be used here.
For $i=\pm 1$ we have a finite generating set $X_i$ of $A_i$ which contains a finite generating set ${\Cal C}$ of $C$.
Then ${\Cal G}=X_1\cup X_{-1}$ is a finite generating set for $G$. 
In section 1 we constructed a language $L_i$ over $X_i$ such that $T_i=\overline {L_i}$ is a right transversal for $C$ in $A_i$.
We also denoted by $\hat T$ the Bass-Serre tree corresponding to presentation $(7)$.
There are two distinguished vertices $d_1=A_1$ and $d_{-1}=A_{-1}$ in $\hat T$. Every positive edge of $\hat T$ is labelled by an element of $T_i$.
Each vertex $v$ of $\hat T$ has the form $s_vA_i$ where $s_v$ is the label of a reduced path in $\hat T$ from $d_{1}$ to $v$.

There is a minimal $H$-invariant subtree $T$ of $\hat T$. Since $H$ is not conjugate to a subgroup of $A_i$, the tree $T$ has at least one edge. We constructed the "fundamental domain" $Y$ for the action of $H$ on $T$ and a finite subtree $Y_1$ of $Y$ which define the algebraic structure of $H$ as the fundamental group of a graph of groups.

By conjugating $H$ we may assume that $(d_1,d_{-1})$ is an edge of $Y_1$. (Notice that a conjugate $H_1$ of $H$ is quasiconvex in $G$ if and only if $H$ is quasiconvex in $G$. Besides, $H_1$ still satisfies the property that $g^{-1}H_1g\cap A_i$ is quasiconvex in $A_i$ for each $g\in G$ and $i=\pm 1$.)

In section 1 we defined a finite graph of groups ${\Bbb B}$ such that $Y_1$ is the maximal subtree of $B$ and $H$ has the structure of the fundamental group of a graph of groups:

$$H=\pi_1({\Bbb B}, Y_1)\eqno (8)$$

We make several immediate observations about $H$.
\proclaim {Lemma 3.2} The group $H$ is word hyperbolic.
\endproclaim
\demo {Proof}
Notice first that each vertex group $B_v$ of ${\Bbb B}$ is word hyperbolic.
Indeed, for $v=s_vA_i\in V{\Bbb B}=VY_1$ we have $B_v=s_vA_is_v^{-1}\cap H\cong A_i\cap s_v^{-1}Hs_v=A_v$.
We know that $A_i\cap s_v^{-1}Hs_v$ is quasiconvex in $A_i$ and thus word hyperbolic. 
Thus we know that $H$ is the fundamental group of a finite graph of groups with virtually cyclic edge groups and word hyperbolic vertex groups. By the results of M.Bestvina and M.Feign [5] and O.Kharlampovich and A.Myasnikov [13], such $H$ is word hyperbolic if and only if it does not contain Baumslag-Solitar subgroups. But $H$ is a subgroup of $G$ which is word hyperbolic and so does not contain Baumslag-Solitar subgroups. Therefore $H$ is word hyperbolic.
\enddemo

We now return to the proof of Theorem B.
For each $v=s_vA_i\in VY_1=VB$ we fix a finite generating set $R_v$ for $A_v$ and a finite generating set $Z_v=s_vR_vs_v^{-1}$ for $B_v=s_vA_vs_v^{-1}$. Recall that each edge $e\in E(B-T)$ is identified with an element $\rho_i(t)\rho_{-i}^{-1}(t)\in H$.
Thus $Z=\underset {v\in VY_1}\to\cup Z_v \bigcup E(B-T)$ is a finite generated set for $H$. Each $e\in EY_1$ represents the trivial element of $H$ so that $Z'=Z\cup EB$ is also a finite generating set for $H$.

Let $h\in H\cap C$. Then $h\in B_{d_1}=A_1\cap H$.
Since $B_{d_1}$ is quasiconvex in $H$ and $A_1$ is quasiconvex in $G$, there is a constant $\lambda_0>0$  such that for any $g\in B_{d_1}$ we have $l_{Z}(g)\le \lambda_0 l_{\Cal G}(g)+\lambda_0$. In particular, $l_{Z}(h)\le \lambda_0 l_{\Cal G}(h)+\lambda_0$. 
We want to establish a similar inequality for the case when $h\not\in C$.

So assume $h\in H-C$. By Proposition B there is a $K'$-quasigeodesic word $W$ with respect to $d_{Z'}$ such that $W$ corresponds to the normal form of $h$ with respect to presentation $(8)$ and $K'$ is a constant independent of $h$.
Thus $W$ has the form

$$W=e_1\dots e_kU_1e_{k+1}\dots  e_{t-1}U_qe_t\dots e_s$$ where 
$$p=e_1,\dots ,e_k,\overline {U_1}, e_{k+1},\dots ,e_{t-1},\overline U_q, e_t,\dots ,e_s$$
is the normal form for $h$ with respect to presentation $(8)$.

Recall that each $\overline {U_i}$ is either an edge of $B-Y_1$ or a nontrivial element of a vertex group of $B$, $e_i\in EY_1$ and $p$ is a loop at the basepoint $d_1$ in the graph of groups ${\Bbb B}$ representing $H$. Recall also that $p$ contains no backtrackings and that no nontrivial element of a vertex group can be "pulled to the left".

Then there is a normal form of $h$ with respect to presentation $(7)$
$$h=v_1\dots v_r$$ satisfying the requirements of Proposition A.

Thus by Proposition A(ii) the syllable length $r$ of $h$ is at least $q_1$ where $q_1$ in the number of those $U_i$ which represent edges of $B-Y_1$.
We also know that for each core element $v_{i_k}\in A_i$ there is a corresponding $U_{j_k}\in B_{v_k}$ such that $l_{Z_{v_k}}(U_{j_k})\le K_0l_{X_i}(\overline {v_{i_k}})$ where $K_0>0$ is some constant independent of $h$.
For each $v_k\in A_i$, $k=1,\dots ,r$ take a $d_{X_i}$-geodesic word $\hat v_k$ representing $v_k$. Put $\hat v=\hat v_1\dots \hat v_r$.

The word $U=U_1\dots U_q$ is obtained from $W$ by deleting some pieces representing identity whose length is bounded by the number of edges in $Y_1$. Thus $U$ is a $K_1$-quasigeodesic with respect to $d_Z$ for some $K_1>0$ independent of $h$. Then 

$$l(U)=m_1+\sum l_{Z_{v_k}}(U_{j_k})\le n + K_0l(\hat v)\le l(\hat v)+K_0l(\hat v)=(K_0+1)l(\hat v)\eqno (9)$$

\proclaim {Calim } We claim that there is a number $N>0$ independent of $h$ such that the following holds.
Suppose $1\le k <n-1$ and $v_k\in A_j-C$, $v_{k+1}\in A_{-j}-C$.
Express $v_k$ as $p_kz_k$ where $z_k\in C$ and $p_k$ is shortest with respect to $d_{X_j}$ in the coset class $v_kC$. Express $v_{k+1}$ as $x_kt_k$ where $x_k\in C$ and $t_k$ is shortest with respect to $d_{X_{-j}}$ in the coset class $Cv_{k+1}$. Then

$$ l_{\Cal C}(z_kx_k)\ge  l_{\Cal C}(z_k)+l_{\Cal C}(x_k)-N\eqno (10)$$ 
\endproclaim  

It is obvious that $(10)$ is satisfied when $C$ is finite. So assume $C$ is infinite. Since $C$ is virtually cyclic and infinite, there is a element $c\in C$ of infinite order and a constant $P>0$ such that for any $c'\in C$ there are integers $n$, $m$ and elements $c_1,c_2\in C$ such that $l_{\Cal C}(c_1), l_{\Cal C}(c_2)\le P$ and $c'=c^nc_1=c_2c^m$.

It is clear that   $l_{\Cal C}(z_k)+l_{\Cal C}(x_k)-l_{\Cal C}(z_kx_k)$ is bounded when at least one of $v_k$, $v_{k+1}$ is not a core element because the lengths of syllables which are not core elements are bounded by 2K (see Proposition A(iii)).
Assume now that they are both core elements.
Thus both $v_k$ and $v_{k+1}$ correspond to the vertex group elements in $U$. Proposition A shows that it is only possible in the following three cases.
\smallskip
{\bf Case 1.} 
There is $U_s=s_vbs_v^{-1}$ and $U_{s+1}=s_was_w^{-1}$ where $v=s_vA_i\in VY_1$, $w=s_vb_1A_{-i}$, $b_1\in T_i$, $b\in A_i$, $a\in A_{-i}$ and
$v_k=v_{i_s}=fbb_1$, $v_{k+1}=v_{i_{s+1}}=af'$ where $f,f'$ have length at most $K$.
\smallskip
{\bf Case 2.} There is $U_s=s_wbs_w^{-1}=s_vb_1bb_1^{-1}s_v^{-1}$ and $U_{s+1}=s_vas_v^{-1}$, where $w=s_wA_{-i}=s_vb_1A_i$, $v=s_vA_i$, $b_1\in T_i$, $b\in A_{-i}\cap A_w$, $a\in A_i\cap A_v$ and
$v_k=v_{i_s}=fb$, $v_{k+1}=v_{i_{s+1}}=b_1af'$ where $f,f'$ have length at most $K$.
\smallskip
We will treat Case 1 and it will be clear that Case 2 is exactly analogous. So assume Case 1 takes place.
Recall that $b\in A_v$, $a\in A_w$ and that the edge group of ${\Bbb B}$ corresponding to the edge $(v,w)$ is $s_vb_1(A_w\cap C)b_1^{-1}s_v^{-1}$.

Recall that $v_k=v_{i_s}=fbb_1$, $v_{k+1}=v_{i_{s+1}}=af_1$, $U_s=s_vbs_v^{-1}$ and $U_{s+1}=s_was_w^{-1}=s_vb_1ab_1^{-1}s_v^{-1}$. Besides we know that
the lengths of $f$ and $f_1$ are bounded by $K$.
We have $v_k=p_kz_k$ where $z_k\in C$ and $p_k$ is shortest with respect to $d_{X_j}$ in the coset class $v_kC$. Likewise, $v_{k+1}=x_kt_k$ where $x_k\in C$ and $t_k$ is shortest with respect to $d_{X_{-j}}$ in the coset class $Cv_{k+1}$.

Recall that $U=U_1\dots U_n$ is a $K_1$-quasigeodesic representative of $h$ in the $d_Z$-metric where $K_1$ does not depend on $h$. In particular $U_sU_{s+1}$ is $K_1$-quasigeodesic. Observe that $U_s$ is a word over $Z_v=s_vR_vs_v^{-1}$ and $U_{s+1}$ is a word over $Z_w=s_vb_1R_wb_1^{-1}s_v^{-1}$. Recall also that $\overline {U_s}=s_vbs_v^{-1}$ and $\overline {U_{s+1}}=s_was_w^{-1}$.
Let $\alpha$ be a $d_{R_w}$-geodesic representative of $a$ and $\beta$ be a $d_{R_v}$-geodesic representative of $b$. 

Express $x_k$ as $x_k=c^nc_a$ and $z_k=c_bc^m$ where $l_{\Cal C}(c_a), l_{\Cal C}(c_b)\le P$. Find a $d_{X_i}$-geodesic representative $u_b$ of $p_kc_b$ and a $d_{X_{-i}}$-geodesic representative $u_a$ of $c_at_k$. Then $v_{k}=fbb_1=\overline {u_bc^m}$ and $v_{k+1}=af_1=\overline {c^nu_a}$. Moreover, since the lengths of $c_a$ and $c_b$ are bounded, Lemma 2.2 implies that there is $\lambda >0$ independent of $h$ such that $u_bc^m$ and $c^nu_a$ are $\lambda$-quasigeodesic in $d_{X_i}$ and $d_{X_{-i}}$ respectively.

We may also assume that $\lambda$ is big enough so that any $R_v$-geodesic word defines a $\lambda$-quasigeodesic in $d_{X_i}$ and any $R_w$-geodesic word defines a $\lambda$-quasigeodesic in $d_{X_{-i}}$.

Thus there is $\epsilon >0$ independent of $h$ such that
\roster
\item "(i)" for any terminal segment $c^{m'}$ of $u_bc^m$ there is a terminal segment $\beta '$ of $\beta$ with $l_{X_i}(\overline {\beta '}b_1c^{-m'})\le \epsilon $ and
\item "(ii)" for any initial segment $c^{n'}$ of $c^nu_a$ there is an initial segment $\alpha '$ of $\alpha$ such that $l_{X_{-i}}(\overline {\alpha '}^{-1}c^{n'})\le \epsilon$.
\endroster
Let $N_1$ be the maximal $l_Z$-length of those elements of $H$ whose ${\Cal G}$-length is at most $2l_{\Cal G}(s_v)+2\epsilon$.
Let $K_2>0$ be such that for any integer $j$ $l_{X_1}(c^j)\ge K_2|j|$ and $l_{X_{-1}}(c^j)\ge K_2|j|$.

Suppose $$|n|+|m|-|n+m| > {{\lambda (K_1N_1+K_1+2\lambda)}\over {K_2}}.$$
Then there is a terminal segment $c^l$ of $u_bc^m$ and an initial segment $c^{-l}$ of $c^nu_a$ such that $$|l|>{{\lambda (K_1N_1+K_1+2\lambda)}\over {2K_2}}.$$ Thus there is a terminal segment $\beta_1$ of $\beta$ and an initial segment $\alpha_1$ of $\alpha$ such that $l_{X_i}(\overline {\beta_1}b_1c^{-l})\le \epsilon $ and $l_{X_{-i}}(\overline {\alpha_1}^{-1}c^{-l})\le \epsilon$.

Therefore $l_{\Cal G}(s_v\overline {\beta_1}s_v^{-1}\cdot s_vb_1\overline {\alpha_1}b_1^{-1}s_v^{-1})\le 2l_{\Cal G}(s_v)+l_{\Cal G}(\overline {\beta_1}b_1\overline {\alpha_1})\le 2l_{\Cal G}(s_v)+2\epsilon +l_{\Cal G}(c^lc^{-l})=2l_{\Cal G}(s_v)+2\epsilon$.
Thus $l(\beta_1\alpha_1)\le K_1l_Z(s_v\overline {\beta_1}s_v^{-1}\cdot s_vb_1\overline {\alpha_1}b_1^{-1}s_v^{-1})+K_1\le K_1N_1+K_1$ since $U$ is a $K_1$-quasigeodesic with respect to $d_Z$.

On the other hand $l(\alpha_1\beta_1)=l(\alpha_1)+l(\beta_1)\ge (1/\lambda)(l_{X_i}(c^l)+\epsilon)-\lambda +(1/\lambda)(l_{X_{-i}}(c^{-l})+\epsilon)-\lambda\ge (2K_2|l|/\lambda)-(2\epsilon/\lambda)> K_1N_1+K_1$ by the choice of $l$.
This gives us a contradiction.
So $\displaystyle{|n|+|m|-|n+m| \le {{\lambda (K_1N_1+K_1+2\lambda )}\over {K_2}}}$ and $(10)$ follows.
\smallskip
{\bf Case 3.} There is $U_s=s_vbs_v^{-1}$, $U_{s+1}=\rho_i(g)\rho_{-i}(g)^{-1}=s_va_ia_{-i}s_w^{-1}$ and $U_{s+2}=s_was_w^{-1}$ such that
$v_k=v_{i_s}=fba_i$ and $v_{k+1}=v_{i_{s+2}}=a_{-i}af'$ where $f, b, a_i\in A_i$, $b_{-i},b,f'\in A_{-i}$ and $f,f'\in \Sigma$.

We have $v_k=p_kz_k$ where $z_k\in C$ and $p_k$ is shortest with respect to $d_{X_i}$ in the coset class $v_kC$. Likewise, $v_{k+1}=x_kt_k$ where $x_k\in C$ and $t_k$ is shortest with respect to $d_{X_{-i}}$ in the coset class $Cv_{k+1}$.

Recall that $U=U_1\dots U_n$ is a $K_1$-quasigeodesic representative of $h$ with respect to $d_{Z}$. In particular $U_sU_{s+1}U_{s+2}$ is $K_1$-quasigeodesic. Observe that $U_s$ is a word over $Z_v=s_vR_vs_v^{-1}$ and $U_{s+2}$ is a word over $Z_w=s_wR_ws_w^{-1}$. Recall also that $\overline {U_s}=s_vbs_v^{-1}$ and $\overline {U_{s+2}}=s_was_w^{-1}$. Let $\beta$ be a $d_{R_v}$-geodesic representative of $b$ and let $\alpha$ be a $d_{R_w}$-geodesic representative of $a$.

Let $x_k=c^nc_a$ and $z_k=c_bc^m$ where $l_{\Cal C}(c_a), l_{\Cal C}(c_b)\le P$.
 Find a $d_{X_i}$-geodesic representative $u_b$ of $p_kc_b$ and a $d_{X_{-i}}$-geodesic representative $u_a$ of $c_at_k$. Then $v_{k}=fba_i=\overline {u_bc^m}$ and $v_{k+1}=a_{-i}^{-1}af_1=\overline {c^nu_a}$. Moreover, there is $\lambda >0$ independent of $h$ such that $u_bc^m$ and $c^nu_a$ are $\lambda$-quasigeodesic in $d_{X_i}$ and $d_{X_{-i}}$ respectively.

We may also assume that $\lambda$ is big enough so that any $R_v$-geodesic word defines a $\lambda$-quasigeodesic in $d_{X_i}$ and any $R_w$-geodesic word defines a $\lambda$-quasigeodesic in $d_{X_{-i}}$.

Thus there is $\epsilon >0$ independent of $h$ such that
\roster
\item "(i)" for any terminal segment $c^{m'}$ of $u_bc^m$ there is a terminal segment $\beta '$ of $\beta$ with $l_{X_i}(\overline {\beta '}a_ic^{-m'})\le \epsilon $ and
\item "(ii)" for any initial segment $c^{n'}$ of $c^nu_a$ there is an initial segment $\alpha '$ of $\alpha$ such that $l_{X_{-i}}(\overline {\alpha '}^{-1}a_{-i}^{-1}c^{n'})\le \epsilon$.
\endroster
Let $N_1$ be the maximal $l_Z$-length of those elements of $H$ whose ${\Cal G}$-length is at most $2l_{\Cal G}(s_v)+2\epsilon$.
Let $K_2>0$ be such that for any integer $j$ we have $l_{X_1}(c^j)\ge K_2|j|$ and $l_{X_{-1}}(c^j)\ge K_2|j|$.
Suppose $$|n|+|m|-|n+m| > {{\lambda (K_1N_1+K_1+2\lambda)}\over {K_2}}.$$
Then there is a terminal segment $c^l$ of $u_bc^m$ and an initial segment $c^{-l}$ of $c^nu_a$ such that $$|l|>{{\lambda (K_1N_1+K_1+2\lambda)}\over {2K_2} } .$$ Thus there is a terminal segment $\beta_1$ of $\beta$ and an initial segment $\alpha_1$ of $\alpha$ such that $l_{X_i}(\overline {\beta_1}a_ic^{-l})\le \epsilon $ and $l_{X_{-i}}(\overline {\alpha_1}^{-1}a_{-i}^{-1}c^{-l})\le \epsilon$.
 
Therefore $l_{\Cal G}(s_v\overline {\beta_1}s_v^{-1}\cdot s_va_ia_{-i}^{-1}s_w^{-1}\cdot s_v\overline {\alpha_1}s_v^{-1})\le 2l_{\Cal G}(s_v)+l_{\Cal G}(\overline {\beta_1}a_ia_{-i}^{-1}\overline {\alpha_1})\le 2l_{\Cal G}(s_v)+2\epsilon +l_{\Cal G}(c^lc^{-l})=2l_{\Cal G}(s_v)+2\epsilon$.
Thus $l(\beta_1)+1+l(\alpha_1)\le K_1l_Z(s_v\overline {\beta_1}s_v^{-1}\cdot s_va_ia_{-i}^{-1}s_w^{-1}\cdot s_w\overline {\alpha_1}bs_w^{-1})+K_1\le K_1N_1+K_1$ since $U$ is a $K_1$-quasigeodesic with respect to $d_Z$.

On the other hand $l(\alpha_1)+1+l(\beta_1)\ge (1\lambda)(l_{X_i}(c^l)+\epsilon)-\lambda +1 +(1\lambda)(l_{X_{-i}}(c^{-l})+\epsilon)-\lambda\ge (2K_2|l|/\lambda)+(2\epsilon/\lambda)-2\lambda +1> K_1N_1+K_1$ by the choice of $l$.
This gives us a contradiction.
So $\displaystyle{|n|+|m|-|n+m| \le {{\lambda (K_1N_1+K_1+2\lambda)}\over {K_2}}}$ and $(10)$ follows.
\smallskip
Thus ve have verified (10) and established the Claim.
\smallskip
Therefore by Lemma 3.1(2) there is a constant $K_3>0$ independent of $h$ such that

$$l_{\Cal G}(h)\ge K_2l(\hat v) -K_2\eqno (11)$$

Recall that $U$ is $K_1$-quasigeodesic with respect to $d_Z$. Thus (9) and (11) imply that 

$$\eqalign{l_{\Cal G}(h)&\ge K_2l(\hat v) -K_2\ge  (K_2/(K_0+1))l(U) -K_2\ge\cr  &(K_2/(K_0+1))K_1l_Z(h)-K_2-(K_2/(K_0+1))K_1}$$
which implies that $H$ is quasiconvex in $G$.

This completes the proof of Theorem B.

\proclaim {Theorem A}
Suppose $G=A_1\ast_C A_{-1}$ is a word hyperbolic group where $C$ is virtually cyclic and $A_1$, $A_{-1}$ have property (Q). Then $G$ has property (Q)
\endproclaim
\demo {Proof}
Let $H$ be a finitely generated subgroup of $G$. Since $C$ is virtually cyclic, for any $g\in G$, $i=\pm 1$ the group $g^{-1}Hg\cap A_i$ is finitely generated. Since $A_i$ has property (Q), the subgroup $g^{-1}Hg\cap A_i$ is quasiconvex in $A_i$. Therefore by Theorem A the subgroup $H$ is quasiconvex in $G$.
\enddemo

\proclaim {Corollary 3.3 (c.f. Corollary 1 from the Introduction)} Let  $G=A_1\ast_C A_{-1}$ where the the groups $A_1, A_{-1}$ belong to the class (Q), $A_1$ is torsion-free and $C$ is a maximal cyclic subgroup of $A_1$. Then $G$ has property (Q).
\endproclaim
\demo {Proof}
By the results of [8], the subgroup $C$ is malnormal in $A_1$. Moreover, since $C$ is cyclic, it is quasiconvex in $A_1$ and in $A_{-1}$. Therefore by the combination theorem for word hyperbolic groups (see [5], [13], [17]) the group $G$ is word hyperbolic. Theorem A implies that $G$ has property (Q).
\enddemo

\proclaim {Corollary 3.4 (c.f. Corollary 2 from the Introduction)} Let $G=A_1\ast_C A_{-1}$ where $C$ is finite and $A_1$, $A_{-1}$ have property (Q). Then $G$ has property (Q).
\endproclaim 
\demo {Proof}
It is easy to show (see [9]) that $G$ is word hyperbolic.
The group $C$ is finite and therefore it is virtually cyclic. Theorem A implies that $G$ has property (Q).
\enddemo

\proclaim {Corollary 3.5 (c.f. Corollary 3 from the Introduction)} Let $G$ be a torsion-free hyperbolic group with property (Q) (e.g. finitely generated free group, hyperbolic surface group etc).
Let $G^{\Bbb Q}$ be the tensor ${\Bbb Q}$-completion of $G$ where ${\Bbb Q}$ is the ring of rational numbers.
Then
\roster
\item  $G^{\Bbb Q}$ is a locally (Q)-group that is any finitely generated subgroup of $G^{\Bbb Q}$ is word hyperbolic and has property (Q);
\item $G^{\Bbb Q}$ has the Howson property that is the intersection of any two finitely generated subgroups of $G$ is finitely generated
\item if $H_1$ and $H_2$ are infinite commensurable subgroups of $G$, that is the intersection $H=H_1\cap H_2$ has finite index in both $H_1$ and $H_2$, then $H$ has finite index in their join $E=gp(H_1\cup H_2)$.
\endroster
\endproclaim
\demo{Proof}

All maximal abelian subgroups of a torsion free word hyperbolic group $G$ are infinite cyclic and malnormal [13].
Therefore, by the results of A.Myasnikov and V.Remeslennikov [15], there is a sequence of groups $G=G_0\le G_1\le\dots \le G_n\le\dots$ such that
\roster
\item "(a)" $G^{\Bbb Q}=\underset {n\ge 0}\to\cup G_n$
\item "(b)" for each $n\ge 0$ we have $G_{n+1}=G_n\underset{g_n=x^{k_n}}\to\ast <x>$ where $g_n\in G_n$ is nontrivial and not a proper power.
\endroster

We claim that each $G_n$ is torsion-free, word hyperbolic and has property (Q). Indeed, it is true for $G=G_0$.
Suppose the claim has been proven for $G_n$, $n\ge 0$.
Then the cyclic subgroup $<g_n>$ is malnormal in $G_n$ since $G_n$ is torsion-free word hyperbolic and $g_n$ is not a proper power (see [13]). The infinite cyclic group $<x>$ is word hyperbolic and has property (Q). The group $G_n$ has property (Q) by inductive hypothesis. Thus by Corollary 3.3 the group $G_{n+1}$ is word hyperbolic and has property (Q). Obviously, $G_{n+1}$ is torsion-free. This concludes the inductive step and the claim is proved.

Any finitely generated subgroup $H$ of $G^{\Bbb Q}$ is contained in some $G_n$, $n\ge 0$ and therefore $H$ has property $Q$. This proves (1). Moreover, if $H_1$ and $H_2$ are finitely generated subgroups of $G^{\Bbb Q}$ then there is $n\ge 0$ such that $H_1, H_2\le G_n$. The group $G_n$ has the Howson property since it belongs to class (Q). Thus $H_1\cap H_2$ is finitely generated. Moreover, if the subgroups $H_1$ and $H_2$ of $G_n$ are commensurable then by the result of [14] their intersection has finite index in their join. This proves (2) and (3).

\enddemo
\proclaim {Corollary 3.6 (c.f. Corollary 5 from the Introduction)}
Suppose $G$ is a one-relator group $G=<x_1,\dots ,x_k,y_1,\dots ,y_s | vu=1>$ where $v$ is a nontrivial freely reduced word in $x_1,\dots ,x_k$ and $u$ is a nontrivial freely reduced word in $y_1,\dots y_s$ which is not a proper power in the free group $F(y_1,\dots ,y_s)$.
Then $G$ has property (Q).
\endproclaim 
\demo {Proof}
Finitely generated free groups $F(x_1,\dots ,x_k)$ and $F(y_1,\dots ,y_s)$ have property (Q) (see [19]). The group $G$ is an amalgamated free product $G=F(x_1,\dots ,x_k)\underset {u^{-1}=v}\to\ast F(y_1,\dots ,y_s)$. Therefore by Corollary 3.3 the group $G$ has property (Q).
\enddemo


\Refs

\ref\no 1
\by J.Alonso, T.Brady, D.Cooper, V.Ferlini, M.Lustig,.M.Mihalik, M.Shapiro and H.Short
\paper Notes on hyperbolic groups
\inbook Group theory from a geometric viewpoint
\bookinfo Proc. ICTP. Trieste
\publ World Scientific\publaddr Singapore
\yr 1991\pages 3--63
\endref


\ref\no 2
\by B.Baumslag
\paper Intersections of Finitely Generated Subgroups in Free Products
\jour J. of London Math. Soc.
\vol 41\yr 1966\pages 673--679
\endref

\ref\no 3
\by G.Baumslag, S.Gersten, M.Shapiro and H.Short
\paper Automatic groups and amalgams
\jour J. of Pure and Appl. Algebra
\vol 76\yr 1991\pages 229--316
\endref

\ref \no 4
\by R.Burns
\paper On the finitely generated subgroups of an amalgamated product of two groups
\jour Trans. Am. Math. Soc.\vol 169\yr 1972\pages 293--306
\endref

\ref \no 5
\by M.Bestvina and M.Feign
\paper The Combination Theorem for Negatively Curved Groups
\jour J. of Diff. Geom.
\vol 35\yr 1992\pages 85-101
\endref

\ref\no 6
\by D.Cohen
\paper Finitely generated subgroups of amalgamated free products and HNN groups
\jour J. of Austr. Math. Soc. Series A\vol 22\issue 3\yr 1976\pages 274--281
\endref

\ref\no 7
\by D.B.A.Epstein,J.W.Cannon,D.F.Holt,S.V.F.Levy, M.S.Paterson and 
W.P.Thurston
\book Word Processing in Groups
\publ Jones and Bartlett, MA \yr 1992
\endref 


\ref\no 8
\by E. Ghys and P. de la Harpe (editors)
\book Sur les groupes hyperboliques d'apr\'es Mikhael Gromov
\publ Birkh\"auser, Progress in Mathematics series, vol. 83
\yr 1990
\endref


\ref\no 9
\by S.Gersten and H.Short
\paper Rational subgroups of biautomatic groups
\jour Ann. Math.
\vol 134\yr 1991\pages 125--158
\endref


\ref \no 10
\by  M.Gromov 
\book Hyperbolic Groups
\bookinfo in 'Essays in group theory', edited by S.M.Gersten
\publ MSRI Publ. 8, Springer \yr 1987 \pages 75-263
\endref 

\ref\no 11
\by I. Kapovich
\paper On a theorem of B.Baumslag
\jour Comm. in Algebra, submitted
\endref

\ref\no 12
\by I.Kapovich
\paper Amalgamated products and the Howson property
\jour Can. Math. Bull., submitted
\endref

\ref\no 13
\by O.Kharlampovich and A.Myasnikov
\paper Hyperbolic groups and free constructions
\jour Trans. Am. Math. Soc., to appear
\endref

\ref\no 14
\by I.Kapovich and H.Short
\paper Greenberg's theorem for quasiconvex subgroups of word hyperbolic groups
\jour Can. J. of Math., submitted
\endref

\ref\no 15
\by A.Myasnikov and V.Remeslennikov
\paper Exponential groups II: extensions of centralizers and and tensor completions of CSA-groups
\jour Intern. J. Alg. Comput. (to appear)
\endref


\ref\no 16
\by S.MacLane
\paper A proof of a subgroup theorem for free products
\jour Matematika\vol 5\yr 1958\pages 13--19
\endref

\ref\no 17
\by P.Papasoglu
\book Geometric methods in group theory
\bookinfo PhD thesis, Columbia University
\yr 1993
\endref 


\ref\no 18
\by C.Pittet
\paper Surface groups and quasiconvexity
\inbook Geometric group Theory (Sussex, 1991) vol. 1
\bookinfo London Math. Soc. Lecture Ser., 181
\publ Cambridge Univ. Press\yr 1993\pages 169--175
\endref

\ref\no 19
\by H.Short
\paper Quasiconvexity and a Theorem of Howson's
\inbook Group theory from a geometric viewpoint
\bookinfo Proc. ICTP. Trieste
\publ World Scientific\publaddr Singapore
\yr 1991
\endref

\ref\no 20
\by G.A.Swarup
\paper Geometric finiteness and rationality
\jour J. of Pure and Appl. Algebra\vol 86\yr 1993
\pages 327--333
\endref

\endRefs
\enddocument